\documentclass[11pt]{article}
\date{~}
\usepackage{amsmath,amsthm,amsfonts,amssymb}
\usepackage[all]{xy}
\expandafter\edef\csname :RestoreCatcodes\endcsname{%
   \catcode`\noexpand :=\the\catcode`:%
   \catcode`\noexpand @=\the\catcode`@%
   \catcode`\noexpand /=\the\catcode`/%
   \catcode`\noexpand &=\the\catcode`&%
   \catcode`\noexpand \^^M=\the\catcode`\^^M%
   \catcode`\noexpand \^^I=\the\catcode`\^^I%
   \expandafter\let\csname:RestoreCatcodes\endcsname=\noexpand\undefined}
\catcode`\:11  \catcode`\@11
   \let\:wlog\wlog \def\wlog#1{}
   \def\:wrn#1#2{\immediate\write\sixt@@n{--DraTeX warning--
      \ifcase #1
    DraTex.sty already loaded
\or \string\Draw\space within \string\Draw
\or Changing definition of \string#2
\or No intersection points: #2
\or Improper rotation of axes: #2
\or (#2) in \string\DSeg\space is a point
\fi}}
\def\:err#1#2{\errmessage{--DraTeX error-- \ifcase #1
     \string#2\space meaningless in three dimensions
\or  \string#2\space meaningless in two dimensions
\or  No \string\MarkLoc(#2)
\or  \string#2 in three dimensions
\or  Too many parameters in definition
\or  \string\MoveFToOval(#2)?
\fi}}
   \ifx\:Xunits\:undefined \else \:wrn0{} \fi
   \catcode`\ 9  \catcode`\^^M9 \catcode`\^^I9
      \newdimen\:LBorder \newdimen\:RBorder\chardef\:eight=8
\mathchardef\:cccvx=360
\newdimen\:mp    \:mp   0.1\p@
\newdimen\:mmp   \:mmp  0.01\p@
\newdimen\:mmmp  \:mmmp 0.001\p@
\newdimen\:XC    \:XC   90\p@
\newdimen\:CVXXX \:CVXXX180\p@
\newdimen\:CCCVX \:CCCVX\:cccvx\p@ \newdimen\:TeXLoc
\newbox\:box\newif\if:IIID  \newdimen\:Z   \newdimen\:Zunits
\newdimen\:Ex   \newdimen\:Ey  \newdimen\:Ez

\def\:AbsVal#1{ \ifdim#1<\z@-\fi #1 }
\def\:abs#1{\ifdim #1<\z@ #1-#1 \fi}
\def\:AbsDif#1#2#3{  #1#2   \advance#1  -#3
   \ifdim #1<\z@ #1-#1 \fi}
\def\:diff#1#2#3{ #1#2  \advance#1 -#3 }
\def\:average#1#2#3{
   #1#2  \advance#1  #3   \divide#1 \tw@}\def\:Opt#1#2#3#4{
   \def\:temp{
      \ifx      \:next\ifnum \def\:next{#3#1#4#2}
      \else\ifx \:next#1     \def\:next{#3}
      \else                  \def\:next{#3#1#4#2}\fi\fi \:next}
   \futurelet\:next\:temp}\def\Define#1{\:multid#1
   \:Opt(){\:Define#1}0}

\def\:DraCatCodes{\catcode`\ 9   \catcode`\^^M9
   \catcode`\^^I9  \catcode`\&13  \catcode`\~13 }

\def\:Define#1(#2){\begingroup  \:DraCatCodes  \::Define#1(#2)}

\def\::Define#1(#2)#3{\endgroup
   \let\:NextDefine\NextDefine
   \let\NextDefine\relax
   \ifcase#2\relax
      \def#1{#3}\or
      \:TxtPar\def#1(##1){#3}\or
         \:TxtPar\def#1(##1,##2){#3}\or
   \:TxtPar\def#1(##1,##2,##3){#3}\or
   \:TxtPar\def#1(##1,##2,##3,##4){#3}\or
   \:TxtPar\def#1(##1,##2,##3,##4,##5){#3}\or
   \:TxtPar\def#1(##1,##2,##3,##4,##5,##6){#3}\or
   \:TxtPar\def#1(##1,##2,##3,##4,##5,##6,##7){#3}\or
   \:TxtPar\def#1(##1,##2,##3,##4,##5,##6,##7,##8){#3}\or
      \:TxtPar\def#1(##1,##2,##3,##4,##5,##6,##7,##8,##9){#3}\or
      \:err4{}\fi      \let\:TxtPar\relax  \:NextDefine}

\let\NextDefine\relax\let\:TxtPar\relax
\def\WarningOn{\def\:multid##1{
   \ifx ##1\:undefined \else \:wrn2##1\fi}}
\def\:gobble#1{}
\def\WarningOff{\let\:multid\:gobble}     \WarningOff
\Define\Indirect{\futurelet\:next\:Indirect}

\Define\:Indirect{\:theDoReg
   \ifx \:next<
      \def\:temp{\let\DoReg\:DoReg}
      \def\:next<##1>{\expandafter\:temp\csname :<##1>\endcsname}
   \else
      \def\:next##1<##2>{
         \expandafter\ifx \csname :<##2> \endcsname \relax
               \def\:next{##1}     \fi
\:indrwrn\Define     \:indrwrn\Object
\:indrwrn\Table      \:indrwrn\IntVar     \:indrwrn\DecVar
         \def\:temp{\let\DoReg\:DoReg##1}
         \expandafter\:temp \csname :<##2> \endcsname}
   \fi      \:next}  \def\:indrwrn#1{    \def\:temp{#1}
   \ifx \:next\:temp \def\:wrn##1##2{\let\:wrn\::wrn} \fi}
\let\::wrn\:wrn
\Define\:Hline{
   \setbox\:box\hbox{\vrule height0.5\:thickness
                    depth 0.5\:thickness width\:x}
   {\:d\:X \advance\:d \wd\:box
 \advance\:X -\:TeXLoc   \global\:TeXLoc\:d
 \vrule width\:X depth\z@ height\z@
 \raise \:Y \box\:box} }\Define\:Vline{
   \setbox\:box\hbox{\vrule width\:thickness
                 \ifdim \:y>\z@ height\:y  depth\z@
                 \else          height\z@  depth-\:y \fi}
   \advance\:X  -0.5\:thickness
   {\:d\:X \advance\:d \wd\:box
 \advance\:X -\:TeXLoc   \global\:TeXLoc\:d
 \vrule width\:X depth\z@ height\z@
 \raise \:Y \box\:box} }\Define\:MvTo(2){\:X#1\:Xunits \:Y#2\:Yunits}
\Define\:Mv(2){\advance\:X  #1\:Xunits
               \advance\:Y  #2\:Yunits}
\def\:DLn(#1,#2,{\:MvTo(#1,#2) \:LnTo(}\Define\:LnTo(2){
   \:x#1\:Xunits \advance\:x  -\:X
   \:y#2\:Yunits \advance\:y  -\:Y
   \:Ln(\:x\du,\:y\du) }\Define\:Ln(2){
   \:x#1\:Xunits \:y#2\:Yunits
   { \ifdim \:x<\z@
        \advance\:X \:x   \:x-\:x
        \advance\:Y \:y   \:y-\:y
     \fi
     \:yy\:AbsVal\:y
     \:dd\:yy \advance\:dd \:x
     \ifdim \:dd>\:mmmp
        \ifdim \:x>\:yy
           { \ifdim\:X<\:LBorder \global\:LBorder\:X\fi
\advance \:X  \:x
\ifdim \:X>\:RBorder \global\:RBorder\:X\fi }
           \let\:Yunitsy\:x  \let\:Xunitsx\:y  \let\:temp\:Hline
           \:dd0.6\:yy    \:divide\:dd\:x
        \else
           { \advance \:X  -0.5\:thickness
\ifdim \:X<\:LBorder \global\:LBorder\:X \fi
\advance \:X  \:thickness
\advance \:X  \:x
\ifdim \:X>\:RBorder \global\:RBorder\:X\fi }
           \let\:Yunitsy\:y  \let\:Xunitsx\:x  \let\:temp\:Vline
           \:dd0.6\:x
\ifdim \:x>\:mmp    \:divide\:dd\:yy   \fi
        \fi
        \advance\:dd  0.4\p@
\:ragged\:Cons\:dd\:ragged
        \:HVLn
      \fi }
   \advance\:X  \:x   \advance\:Y  \:y}\Define\:HVLn{
   \:xx\:AbsVal\:Xunitsx  \:divide\:xx\:ragged
   \advance\:xx  0.99\p@   \:K\:InCons\:xx \relax
   \ifnum \:K>\z@
      \divide\:Xunitsx \:K  \advance\:K \@ne
      \divide\:Yunitsy \:K
   \else \:K\@ne  \fi  \:NextLn}

\Define\:NextLn{
   \ifnum\:K=\z@  \let\:NextLn\relax
   \else  { \:temp }  \advance\:K  \m@ne
      \advance\:X  \:x   \advance\:Y  \:y
   \fi   \:NextLn}\newdimen\:ragged

\Define\Ragged(1){    \:ragged#1\p@  \:ragged0.1\:ragged }
\Ragged(7.5)
\Define\PaintUnderCurve(4){{
   \:Z\:Y   \def\:next{\Curve(#1,#2,#3,#4)}
   \MoveToLoc(#1)  \:d\:X   \MoveToLoc(#4)
   \advance\:d -\:X
   \ifdim \:AbsVal\:d<\:mmp  \def\:next{}
   \else
      \def\:CrvLnTo(##1,##2){
         \:x \:X   \:y  \:Y    \:X\:DJ  \:Y\:yyyy
         \:xx\:X   \:dddd \:Y  \:X\:x  \:Y\:y
         { \advance\:Y  \:dddd   \divide\:Y \tw@
           \advance\:Z -\:Y
           \advance\:Y   0.5\:Z
           \:dddd  \:AbsVal \:Z  \:d\z@
           \def\:CrvLnTo{\:LnTo}
           \:yy\:Y  \:dd\:dddd  \:ddd\:dddd
           \::paint  }}
    \fi \:next       }}
\Define\DoCurve(1){ \let\:StopCurve\:SlowCurve
  \def\:CMv(##1){  \:x\:X \:y\:Y   \:MvTo(##1)
      \advance\:x -\:X   \advance\:y -\:Y
      \:xxx \:x    \:yyy\:y}
   \:DoCurve{\Curve(#1)}
   \let\:StopCurve\:FastCurve}

\def\:DoCurve#1(#2)#3{{\XSaveUnits
   \def\:next{#1}    \:MvTo(#2,#2)
   \:x\:AbsVal\:X  \:y\:Y  \:ddd\z@  \:length
   \:Z\:d   \:divide\:Z{1.41421\p@}
   \edef\:tempa{\the\:DoDist}   \global\:DoDist\z@
   \def\:CrvLnTo(##1){ \MarkLoc(1^)    \:CMv(##1)
      { \MarkLoc(2^)   \:ddd\z@    \:length
\:dd  \:DoDist  \global\advance\:DoDist  \:d
\:ddd \:DoDist  \:divide\:ddd\:Z
\DoReg\:InCons\:ddd  \:Z\DoReg\:Z

      \ifdim \:Z>\:dd
         \advance\:Z -\:DoDist
\advance\:dd -\:DoDist
\:divide\:Z\:dd
\advance\:X \:Cons\:Z\:xxx
\advance\:Y \:Cons\:Z\:yyy   \:DoRot
         \def\:CrvLnTo{\:LnTo}
         \def\:OvalLn{\:Ln}  \XRecallUnits    #3 \fi}}
   \:next    \xdef\:DoDim{\:Cons\:DoDist}
   \global\:DoDist\:tempa     }
   \let\DoDim\:DoDim}

\newdimen\:DoDist\def\:DoRot{ \DSeg\RotateTo(1^,2^) }
\def\DoLine(#1,#2)(#3)#4{
   \MarkLoc($1)  \Move(#1,#2)
   \def\:next{   { \MarkLoc($2)
     \DSeg\RotateTo($1,$2)   \let\:DoRot\relax
     \edef\:RecallRagged{\the\:ragged} \MoveTo(#3,#3)
     \:x\:AbsVal\:X  \:y\:Y  \:ddd\z@  \:length
     \:ragged\:d   \divide\:ragged \tw@
     \DoCurve($1,$1,$2,$2)(#3)
        {\:ragged\:RecallRagged #4}  }
   \let\DoDim\:DoDim}  \:next }
\def\Table#1{\begingroup  \:DraCatCodes   \:multid#1
   \:DefineData#1}

\def\:DefineData#1#2{\endgroup
   \let\:temp~  \def~{\noexpand~}
   \edef#1{\noexpand\:DoPoly\expandafter\noexpand\csname :\string#1\endcsname}
   \expandafter\edef  \csname :\string#1\endcsname
       ##1{\noexpand\ifcase##1(#2)\noexpand\fi}
   \let~\:temp   \:DoNextPoly   \:DoNextPoly}

\def\:OR{\let\:or\or}  \:OR \catcode`\&13  \def&{)\noexpand\:or(}

\def\:TableData#1#2#3{\endgroup   \Table\:temp{#3}
   \:K\z@   \:J\z@    \def\:tempa(##1){\advance\:J \@ne }
   \:temp(0,999){\:tempa}    \let\:tempa&      \def\:temp{\def#1}
   \def&##1&{
      \ifnum  \:K<\:J
         \advance\:K \@ne
         \ifnum  \:K=\@ne   \def#1{#2(##1)}
         \else
            \:IIIexpandafter\:temp\expandafter{
               #1 & #2(##1) }
         \fi
      \else  \let\:next\relax \fi
      \:next}
   \let\:next&     &#3&&   \let&\:tempa }

\catcode`\&4  \def\:DoPoly#1(#2,#3)#4{
   \expandafter\let \csname :Back\the\:level\endcsname\:or
\expandafter\edef\csname :DoVars\the\:level\endcsname{
   \:DoB\the\:DoB}
\advance\:level  \@ne
   \:DoB#3  \advance\:DoB -#2
   \def\:PolyOr(##1){
      \ifnum  \:DoB=\z@  \:OR
      \else   #4(##1)   \advance\:DoB \m@ne    \fi}
   \:OR
   \def\:temp{\let\:or\:PolyOr #4}
   \:IIIexpandafter\:temp#1{#2}
   \advance\:level  \m@ne
\csname :DoVars\the\:level\endcsname
\def\:temp{\let\:or}
\expandafter\:temp\csname :Back\the\:level\endcsname }
\Define\PaintRect(2){\def\:next{{   \MarkLoc(^)
   \MoveToLoc(^) \Move( #1,0) \MarkLoc(^1) \Move(0,#2)
   \MarkLoc(^2)  \Move(-#1,0) \MarkLoc(^3)
   \PaintQuad(^,^1,^2,^3)        }}\:next}

\Define\PaintRectAt(4){\def\:next{{   \MoveTo(#1,#2)
   \:K#3   \advance\:K -#1
   \:DoB#4 \advance\:DoB -#2  \PaintRect(\:K,\:DoB)}}
   \:next}\def\::paint{
   \ifdim \:d<\:ragged      \advance\:xx -\:X
      \:yyy\:Y \:xxx\:dddd
      \advance\:Y \:yy  \divide\:Y \tw@
      \:average\:dddd\:dd\:ddd
      \def\:next{\:brush(\:xx,\z@)\:Y\:yyy\:dddd\:xxx}
   \else  \divide\:d \tw@
      \:average\:x\:X\:xx
      \:average\:y\:Y\:yy
      \:average\:dddd\:dd\:ddd
   \fi   \:next}\Define\:paint{{   \:AbsDif\:d\:xx\:x
   \ifdim \:d<\:mmp  \let\::paint\relax
   \else
      \ifdim  \:y >\:yyy   \:dd\:y   \:y \:yyy   \:yyy \:dd  \fi
      \ifdim  \:yy>\:yyyy  \:dd\:yy  \:yy\:yyyy  \:yyyy\:dd  \fi
      \:AbsDif\:dd\:yyyy\:yyy      \:AbsDif\:ddd\:yy\:y
      \ifdim \:dd<\:ddd  \:dd\:ddd  \fi
      \ifdim \:d >\:dd   \:d\:dd  \fi
      \advance\:d \:d      \:diff\:dd\:y\:yyy
      \:diff\:ddd\:yy\:yyyy
      \advance\:y    -0.5\:dd    \:abs\:dd
      \advance\:yy   -0.5\:ddd   \:abs\:ddd
      \:X\:x  \:Y\:y
      \:average\:dddd\:dd\:ddd
      \def\:next{ \:lpaint \:rpaint }
   \fi
   \::paint }}

\Define\:lpaint{ { \:xx\:x  \:yy \:y  \:ddd\:dddd  \::paint} }
\Define\:rpaint{ { \:X \:x  \:Y  \:y  \:dd \:dddd  \::paint} }
\Define\PaintQuad(4){\def\:next{{\Units(1pt,1pt)
   \MoveToLoc(#1)  \:x  \:X  \:y  \:Y
   \MoveToLoc(#2)  \:xx \:X  \:yy \:Y
   \MoveToLoc(#3)  \:xxx\:X  \:yyy\:Y
   \MoveToLoc(#4)  \:xxxx  \:X  \:yyyy \:Y

   \:paintQuad  }}\:next}

\def\:paintQuad{{
   \:SetVal\:a\:x\:y\:xx\:yy\:xxxx\:yyyy
\:SetVal\:b\:xx\:yy\:xxx\:yyy\:x\:y
\:SetVal\:c\:xxx\:yyy\:xx\:yy\:xxxx\:yyyy
\:SetVal\:cc\:xxxx\:yyyy\:xxx\:yyy\:x\:y
\def\:A{\:a} \def\:B{\:b} \def\:C{\:c} \def\:D{\:cc}
\:sort\:B\:A
\:sort\:C\:B  \:sort\:B\:A
\:sort\:D\:C  \:sort\:C\:B  \:sort\:B\:A
\let\:temp\relax
\:IsTriang\:A\:B
\:IsTriang\:B\:C
\:IsTriang\:C\:D
\:temp
   \:Quad\:A\:temp>   \:xxxx\:xx \:yyyy\:yy \:Z\:d
\:Quad\:D\:next<
   \:PrePaint(\:xx,\:yy,\:d,\:xxxx,\:yyyy,\:Z)
   \:temp  \:next }}
\Define\:PrePaint(6){
   \:x#1 \:y#2 \:yyy#3 \:xx#4 \:yy#5 \:yyyy#6
   \:paint }\def\:Quad#1#2#3{
         \:GetVal#1\:x\:y0
   \:GetVal#1\:xx\:yy1
   \:GetVal#1\:xxx\:yyy2
   \ifdim \:xx#3\:xxx
      \:ddd\:xx   \:xx\:xxx   \:xxx\:ddd
      \:ddd\:yy   \:yy\:yyy   \:yyy\:ddd
   \fi
          \def#2{}
   \:diff\:dd\:xxx\:xx
   \ifdim \ifdim\:AbsVal\:dd<\:mp  \:yy=\:yyy   \else  \z@>\z@  \fi

      \:d\:yy
   \else
      \ifdim \:AbsVal\:dd>\:mp
         \:diff\:dd\:xxx\:x   \:diff\:ddd\:yyy\:y
         \:divide\:ddd\:dd    \:diff\:dd\:xx\:xxx
         \:ddd\:Cons\:ddd\:dd   \advance\:ddd \:yyy
         \:d\:ddd
      \else \:d\:yyy \fi
      \edef#2{  \noexpand\:PrePaint
         (\the\:x ,\the\:y ,\the\:y,
         \the\:xx,\the\:yy,\the\:d)    }
   \fi}\def\:SetVal#1#2#3#4#5#6#7{
   \edef#1{(\the#2,\the#3,\the#4,\the#5,\the#6,\the#7)}}

\def\:sort#1#2{
   \ifdim \:IIIexpandafter\:field#1 <
          \:IIIexpandafter\:field#2
      \let\:temp#1  \let#1#2  \let#2\:temp
   \fi  }

\Define\:field(6){#1}

\def\:GetVal#1#2#3{
   \:IIIexpandafter\::GetVal #1#2#3}

\def\::GetVal(#1,#2,#3,#4,#5,#6)#7#8#9{
      \ifcase #9 #7#1   #8#2\or #7#3   #8#4\or #7#5   #8#6 \fi}
\def\:IsTriang#1#2{
   \ifdim \:IIIexpandafter\:field#1 =
          \:IIIexpandafter\:field#2
      \ifdim \:IIIexpandafter\:fieldB#1 =
             \:IIIexpandafter\:fieldB#2
         \def\:temp{ \:FixTria }
   \fi \fi  }

\def\:FixTria{
   \edef\:temp{\:IIIexpandafter\:FrsII\:B}
   \ifdim \:IIIexpandafter\:field\:A =
          \:IIIexpandafter\:field\:B
      \ifdim \:IIIexpandafter\:fieldB\:A =
             \:IIIexpandafter\:fieldB\:B
          \edef\:temp{\:IIIexpandafter\:FrsII\:C}
   \fi\fi
   \edef\:A{\:IIIexpandafter\:FrsII\:A}
   \edef\:D{\:IIIexpandafter\:FrsII\:D}
   \edef\:temp{
      \def\noexpand\:a{(\:A,\:temp,\:D)}
      \def\noexpand\:b{(\:temp,\:A,\:D)}
      \def\noexpand\:c{\noexpand\:b}
      \def\noexpand\:cc{(\:D,\:A,\:temp)}}
   \:temp
   \def\:A{\:a}  \def\:B{\:b}  \def\:C{\:c}  \def\:D{\:cc}  }

\Define\:fieldB(6){#2}
\Define\:FrsII(6){#1,#2}
\def\:IIIexpandafter{\expandafter\expandafter\expandafter}

\Define\DrawRect(2){ \Line( #1,0) \Line(0, #2)
                     \Line(-#1,0) \Line(0,-#2)}

\Define\DrawRectAt(4){{
   \MoveTo(#1,#2) \LineTo(#3,#2) \LineTo(#3,#4)
   \LineTo(#1,#4) \LineTo(#1,#2)}}
\def\du#1{ \ifx#1\:Xunits   \else\ifx#1\:Yunits
      \else\ifx#1\:Zunits   \else #1
      \fi \fi \fi}\Define\XSaveUnits{
   \expandafter\edef\csname XRecallUnits\the\:level\endcsname{
      \:StoreUnits}
    \advance\:level  \@ne}

\Define\XRecallUnits{
   \advance\:level \m@ne
   \csname XRecallUnits\the\:level \endcsname}

\Define\SaveUnits{
     }

\Define\:StoreUnits{       \:Xunits \the\:Xunits
   \:Yunits \the\:Yunits \:Zunits \the\:Zunits
   \:Xunitsx\the\:Xunitsx \:Xunitsy\the\:Xunitsy
\:Yunitsx\the\:Yunitsx \:Yunitsy\the\:Yunitsy }
\Define\:SearchDir{
   \ifdim \:x<\z@  \:x-\:x  \:y-\:y
   \edef\:tempA{\advance\:ddd  \ifdim \:y<\z@ - \fi\:CVXXX}
\else         \def\:tempA{} \fi
\ifdim \:y<\z@
   \edef\:tempA{\:ddd-\:ddd \advance\:ddd  \:CCCVX \:tempA}
   \:y-\:y
\fi
\ifdim \:y>\:x
   \:ddd\:y \:y\:x \:x\:ddd
   \edef\:tempA{\advance\:ddd  -\:XC \:ddd-\:ddd \:tempA}
\fi
   \:divide\:y\:x  \:d57.29578\:y
   \:ddd\:d       \:sqr\:y  \:K  \@ne
   \Do(1,30){\advance\:K \tw@
      \:d-\:Cons\:y\:d  \:dd\:d  \divide\:dd \:K
      \advance\:ddd \:dd  }
   \advance\:ddd -0.49182\:dd
   \:tempA }\Define\Curve(4){\def\:next{{  \XSaveUnits \Units(1pt,1pt)
   \MoveToLoc(#1) \:DI \:X  \:DK \:Y
   \MoveToLoc(#2) \:ddd \:X  \:Ez \:Y
   \MoveToLoc(#3) \:dd\:X  \:Ey\:Y
   \MoveToLoc(#4) \:DJ\:X  \:yyyy\:Y   \:Curve    }}\:next}
\Define\:FastCurve{   \:AbsDif\:d\:DI\:ddd  \:AbsDif\:dddd\:DK\:Ez
   \advance\:d \:dddd
   \ifnum \:d<\:ragged
      \:AbsDif\:d\:DJ\:dd \:AbsDif\:dddd\:yyyy\:Ey
      \advance\:d \:dddd                     \fi}
\let\:StopCurve\:FastCurve
\Define\:SlowCurve{
   \:AbsDif\:d\:DI\:DJ  \:AbsDif\:dddd\:DK\:yyyy
   \advance\:d \:dddd                     }
\Define\:Curve{   \:StopCurve
   \ifnum \:d<\:ragged   \:X\:DI \:Y\:DK
                    \:CrvLnTo(\:Cons\:DJ,\:Cons\:yyyy)
   \def\:SubCurves{}  \fi
   \:SubCurves}

\def\:CrvLnTo{\:LnTo}\Define\:SubCurves{
   \:average\:yy\:DI\:ddd     \:average\:Ex\:DK\:Ez
   \:average\:ddd\:ddd\:dd    \:average\:Ez\:Ez\:Ey
   \:average\:dd\:dd\:DJ  \:average\:Ey\:Ey\:yyyy
   \:average\:Zunits\:yy\:ddd    \:average\:Vdirection\:Ex\:Ez
   \:average\:ddd\:ddd\:dd    \:average\:Ez\:Ez\:Ey
   \:average\:DL\:Zunits\:ddd  \:average\:xxxx\:Vdirection\:Ez
   { \:ddd \:yy  \:Ez \:Ex   \:dd\:Zunits
     \:Ey\:Vdirection \:DJ\:DL \:yyyy\:xxxx \:Curve }
   \:DI \:DL \:DK \:xxxx               \:Curve   }
\def\MoveToCurve[#1]{
   \Define\:BiSect(3){\MoveToLoc(##1)
      \CSeg[#1]\Move(##1,##2)
   \MarkLoc(##3) }\:MvToCrv}

\Define\:MvToCrv(4){  \:BiSect(#1,#2,:a) \:BiSect(#2,#3,:b)
   \:BiSect(#3,#4,:c) \:BiSect(:a,:b,:A) \:BiSect(:b,:c,:B)
   \:BiSect(:A,:B,:Q)}
\Define\:OvalDir(3){   \:CosSin{#3\p@}
   \:Zunits\p@   \:d#2\:Zunits   \:x\:Cons\:d\:x
   \:Zunits\p@   \:d#1\:Zunits   \:y\:Cons\:d\:y
   \:SearchDir   }\Define\DrawOval   (2){ \DrawOvalArc(#1,#2)(0,\:cccvx) }
\Define\PaintOval  (2){ \PaintOvalArc(#1,#2)(0,\:cccvx) }
\Define\DrawCircle (1){ \DrawOval(#1,#1) }
\Define\PaintCircle(1){ \PaintOval(#1,#1) }
\def\DrawOvalArc(#1,#2)(#3,#4){{
       \:xxxx#4\p@  \advance\:xxxx -#3\p@
       \ifdim \:xxxx=\z@ \else
   \let\:SinOne\:SinB  \:OvalDir(#1,#2,#3)  \:DJ\:ddd
\:OvalDir(#1,#2,#4)  \:diff\:DI\:ddd\:DJ
\ifdim\:DI<\z@ \advance\:DI  \:CCCVX \fi
\ifdim \:xxxx<\:CCCVX \else \:DI\:CCCVX \fi
\:InitOval(#1,#2)  \:CosSin\:DJ
   \:xxxx\:x  \:yyyy\:y   \:xx\:X  \:yy\:Y
   \advance\:X \:Xx\:x  \advance\:X \:Yx\:y
   \advance\:Y \:Xy\:x  \advance\:Y \:Yy\:y
   \let\:Xunits\empty  \let\:Yunits\empty
   \Do(1,\:InCons\:DI){
      \:dd\:X  \:ddd\:Y  \:X\:xx  \:Y\:yy
      \:AdvOv\:xxx\:yyy\:xxxx\:yyyy
      \:X\:dd  \:Y\:ddd
      \advance\:xxx -\:X   \advance\:yyy -\:Y
      \:d\:AbsVal\:xxx
      \advance\:d  \:AbsVal\:yyy
      \ifdim \:d>\:ragged
         \:OvalLn(\:xxx,\:yyy)  \fi  }
   \:OvalDir(#1,#2,#4)    \:CosSin\:ddd
   \advance\:xx \:Xx\:x  \advance\:xx \:Yx\:y
   \advance\:yy \:Xy\:x  \advance\:yy \:Yy\:y
   \advance\:xx -\:X     \advance\:yy -\:Y
   \:OvalLn(\:xx,\:yy) \fi }}

\def\:OvalLn{\:Ln}\def\DoOvalArc(#1)(#2){   \:xx\:X  \:yy\:Y
   \def\:CMv(##1){  \:Mv(\:xxx,\:yyy)
      \:x\:xxx  \:y\:yyy}
   \:DoCurve{             \:X\:xx  \:Y\:yy
       \def\:DoRot{  \let\:Xunits\:XunitsReg
                     \let\:Yunits\:YunitsReg
                     \DSeg\RotateTo(1^,2^)     }
       \let\::OvalLn\:CrvLnTo
\Define\:OvalLn(2){ \:dd\:AbsVal####1
   \advance\:dd \:AbsVal####2 \:divide\:dd\:ragged
   \:J\:InCons\:dd  \advance\:J  \@ne
   \divide####1  \:J  \divide####2  \:J
   \Do(1,\:J){\::OvalLn(####1,####2)}}
       \DrawOvalArc(#1)(#2)}}
\def\NextTable{\begingroup  \:DraCatCodes \:NextTable}

\def\:NextTable#1{\endgroup
  \def\:DoNextPoly{#1\NextTable{}}}
\NextTable{}\def\:AdvOv#1#2#3#4{
   \:d\:CosOne#3 \advance\:d -\:SinOne#4
   #4\:CosOne#4    \advance#4     \:SinOne#3   #3\:d
   \divide#3 \:eight  \divide#4 \:eight
   #1\:X   #2\:Y
   \:d\:Xx#3  \advance\:d \:Yx#4  \advance#1  \:d
   \:d\:Xy#3  \advance\:d \:Yy#4  \advance#2  \:d  }

\def\:CosOne{7.99878}   \def\:SinB{0.13962}
\def\PaintOvalArc(#1,#2)(#3,#4){{ \ifdim #3\p@=#4\p@
                      \let\:next\relax  \else
   \:d\:AbsVal{#1\:Xunits} \advance\:d \:AbsVal{#2\:Yunits}
   \ifdim \:d<3\:ragged   \divide \:d \tw@   \PenSize(\:d)
\:Mv(-0.5\:d\du,0)  \:Ln(\:d\du,0)
   \else    \:InitOval(#1,#2)
      \MarkLoc(o$)   \RotateTo(#3) \MoveFToOval(#1,#2)
      \:Ex\:X  \:Ey\:Y   \edef\:FirstOvDir{\:Cons\:ddd\p@}
      \MoveToLoc(o$) \RotateTo(#4) \MoveFToOval(#1,#2)
      \:Ez \:X  \:Vdirection \:Y   \edef\:LastOvDir{\:Cons\:ddd\p@}
      \MoveToLoc(o$)
      \if:rotated
      \:Zunits\p@     \:Zunits#1\:Zunits
      \:xx\:Cons\:Zunits\:Xunitsx
      \:Zunits\p@     \:Zunits#2\:Zunits
      \:yy\:Cons\:Zunits\:Yunitsx
      \:x\:xx  \:y\:yy  \:ddd\z@  \:length
      \:ddd\:d  \:divide\:xx\:ddd   \:divide\:yy\:ddd
\else \:xx\p@  \:yy\z@   \fi
      \:AbsDif\:d{#3\p@}{#4\p@}
      \ifdim \:d>359\p@ \:Ez-\:Xx\:xx  \advance\:Ez -\:Yx\:yy
\:Vdirection-\:Xy\:xx  \advance\:Vdirection -\:Yy\:yy
\advance\:Ez \:X  \advance\:Vdirection \:Y
\:setpaint\:PaintOvOv<>
      \else
         \:x\:xx \:y\:yy \:SearchDir
\:xxx\:FirstOvDir  \advance\:xxx -\:ddd
\ifdim \:xxx<\z@    \advance\:xxx \:CCCVX  \fi
\:yyy\:LastOvDir   \advance\:yyy -\:ddd
\ifdim \:yyy<\z@    \advance\:yyy \:CCCVX  \fi
\:J\z@
\ifdim \:xxx<\:yyy  \ifdim        \:yyy<\:CVXXX
            \:Pntovln\:FirstOvDir\:FirstOvDir\:LastOvDir
\else \ifdim  \:xxx>\:CVXXX
            \:Pntovln\:FirstOvDir\:FirstOvDir\:LastOvDir
\else \:yyy-\:yyy  \advance\:yyy \:CCCVX
      \ifdim  \:xxx<\:yyy
            \:setpaint\:PntLeftOvOv><  \:PntMovln\:FirstOvDir\:xx
      \else \:FxLx  \:setpaint\:PntLeftOvOv><  \:Pntmovln\:LastOvDir
\fi  \fi  \fi
\else               \ifdim        \:xxx<\:CVXXX
            { \:setpaint\:PaintOvOv<> }  \:FxLx  \:setpaint\:PntLeftOvOv><
\:Usrch  \:PaintMidOvLn\:LastOvDir\:FirstOvDir
\else \ifdim  \:yyy>\:CVXXX
            {  \:FxLx \:setpaint\:PaintOvOv<>  }  \:setpaint\:PntLeftOvOv><
 \:Dsrch  \:PaintMidOvLn\:LastOvDir\:FirstOvDir
\else \:xxx-\:xxx  \advance\:xxx \:CCCVX
      \ifdim  \:yyy<\:xxx
            \:setpaint\:PaintOvOv<>    \:PntMovln\:FirstOvDir\:xx
      \else
            \:FxLx  \:setpaint\:PaintOvOv<>   \:Pntmovln\:LastOvDir
\fi  \fi  \fi  \fi
   \fi \fi \fi}}\def\:setpaint#1#2#3{{\aftergroup#1
\:d\:Xx\:xx   \advance\:d \:Yx\:yy
\ifdim \:d<\z@  \aftergroup#3
\else           \aftergroup#2  \fi}}
\def\:FxLx{\:d\:Ex  \:Ex\:Ez  \:Ez\:d}

\def\:Pntovln#1{
   \let\:SinOne\:SinB     \:CosSin#1
   \:xxxx\:Ex  \:yyyy\:Ey  \:Z\z@
   \:PaintOvLn}

\def\:PntMovln{
   \:Dsrch  \:FxLx   \:xx\:ddd  \:PaintOvLn}

\def\:Pntmovln{
   \:Usrch  \:FxLx  \:xx\:ddd   \:PaintOvLn\:xx}
\def\:PaintMidOvLn#1#2{
   \:FxLx   \:xx\:ddd  \:PaintOvLn#1#2
   \:xxx\:Ez  \:yyy\:Vdirection
   \:ddd\:Xy\:x  \advance\:ddd \:Yy\:y  \advance\:ddd \:Y
   \:PaintSlice}

\def\:PntLeftOvOv{
   \:xx-\:xx  \:yy-\:yy  \:PaintOvOv}

\def\:Dsrch{      \def\:SinOne{-\:SinB}
   \:xx\:x  \:yy\:y   \:SearchDir
   \:x\:xx  \:y\:yy
   \:d\:Ey  \:Ey\:Vdirection  \:Vdirection\:d }

\def\:Usrch{
   \let\:SinOne\:SinB
   \:x\:xx  \:y\:yy   \:SearchDir
   \:x\:xx  \:y\:yy}\def\:PaintOvLn#1#2{
   \:diff\:dd\:Ey\:Vdirection  \:diff\:ddd\:Ex\:Ez
   \ifdim \:AbsVal\:ddd>\:mp
      \:divide\:dd\:ddd      \:d#2
      \advance\:d -#1
      \ifnum \:d<\z@  \advance\:d \:CCCVX  \fi
      \:DoB\:InCons\:d  \let\:next\:PntDo  \:next
   \fi}

\def\:PntDo{
   \ifnum\:DoB=\z@ \let\:next\relax
   \else
      \::AdvOv\:x\:y
      \ifdim \:d>\:ragged
         \:ddd\:xxx  \advance\:ddd -\:Ez
         \:ddd\:Cons\:dd\:ddd
         \advance\:ddd \:Vdirection  \:PaintSlice
      \fi
      \advance\:DoB \m@ne
   \fi  \:next}\Define\:PaintSlice{   \:AbsDif\:dddd\:yyy\:ddd
   \advance\:yyy \:ddd  \divide\:yyy \tw@
   { \advance\:dddd \:Z  \divide\:dddd \tw@
      \:X\:xxxx  \:Y\:yyy
      \:xx\:xxx  \advance\:xxx -\:xxxx
      \:d\z@
      \:yy\:Y  \:dd\:dddd  \:ddd\:dddd
      \::paint     }
   \:xxxx\:xxx  \:yyyy\:yyy  \:Z\:dddd    \:J\z@}
\def\::AdvOv#1#2{  \:AdvOv\:xxx\:yyy#1#2
   \:AbsDif\:d\:xxxx\:xxx       \advance\:J \@ne
   \ifnum \:J=\sixt@@n   \multiply\:d \@cclvi
   \fi }\def\:PaintOvOv#1{   \def\:hdir{#1}
   \:xxx\:Xx\:xx  \advance\:xxx \:Yx\:yy
   \:yyy\:Xy\:xx  \advance\:yyy \:Yy\:yy
   \advance\:xxx \:X   \advance\:yyy \:Y
   \:Z\z@  \:xxxx\:xxx  \:yyyy\:yyy
   \:x\:xx   \:y\:yy   \:DoB\z@   \:J\z@
   \let\:next\:scanOvOv  \:next }
\Define\:scanOvOv{  \advance\:DoB \@ne
   \advance\:J \@ne
   \:AbsDif\:d\:xxx\:Ez
   \ifdim \ifdim      \:xxx\:hdir\:Ez    -
          \else\ifdim \thr@@\:d<\:ragged -
          \else\ifnum \:DoB>358            -
          \fi \fi \fi                     \p@<\z@
      \:ddd\:yyyy  \advance\:ddd -0.5\:dddd
      \:yyy\:yyyy   \advance\:yyy   0.5\:dddd
      \:xxx\:Ez \:PaintSlice
      \let\:next\relax
   \else
     \:d\:xx \advance\:d -\:x
     \:xxx\:yy \advance\:xxx -\:y
     \ifdim -\:Xx\:d\:hdir\:Yx\:xxx
        \let\:SinOne\:SinB
        \::AdvOv\:xx\:yy
        \ifdim  \ifdim       \:d >\:ragged -
                \else \ifnum \:J>\sixt@@n -
                \fi\fi                      \p@<\z@
           \:J\z@  \def\:SinOne{-\:SinB}
           \:AdvOv\:dd\:ddd\:x\:y  \:PaintSlice  \fi
     \else
        \def\:SinOne{-\:SinB}
        \::AdvOv\:x\:y
        \ifdim  \ifdim       \:d >\:ragged -
                \else \ifnum \:J>\sixt@@n -
                \fi\fi                      \p@<\z@
           \:J\z@   \let\:SinOne\:SinB
           \:AdvOv\:dd\:ddd\:xx\:yy  \:PaintSlice  \fi
   \fi \fi
   \:next}\Define\SetBrush{\:Opt[]\:SetBrush{}}

\def\:SetBrush[#1](#2,#3)#4{    \def\:temp{#4}
   \ifx \:temp\empty
      \def\:brush{   \let\:Xunits\empty \let\:Yunits\empty
                     \:thickness\:dddd  \:Ln  }
   \else       \def\:BruShape{#4}
      \:dd#2\:Xunits      \:ddd#3\:Yunits
      \edef\:Grd{ \:dd\the\:dd  \:ddd\the\:ddd }
      \MarkLoc($$)  \def\:temp{#1}
\ifx \:temp\empty  \:X\z@  \:Y\z@
\else  \MoveTo(#1) \fi
\edef\:BrOrg{ \:x\the\:X  \:y\the\:Y }
\MoveToLoc($$)
      \def\:brush(##1,##2){ \::brush }  \fi  }

\SetBrush(,){}\def\::brush{{  \SetBrush(,){}
\let\:Xunits\:XunitsReg   \let\:Yunits\:YunitsReg
\advance\:Y -0.5\:dddd   \:yy\:Y
\advance\:yy \:dddd
\:BrOrg  \:Grd  \advance\:xx  \:X
\ifdim \:xx<\:X  \:d\:X \:X\:xx \:xx\:d  \fi
   \:GridPt\:X\:x\:dd
   \:GridPt\:Y\:y\:ddd            \:x\:X
   \:DoBrush       }}\Define\:DoBrush{
   \ifdim       \:Y>\:yy  \let\:DoBrush\relax
   \else \ifdim \:X>\:xx \advance\:Y \:ddd     \:X\:x
   \else { \:BruShape }  \advance\:X \:dd
   \fi  \fi  \:DoBrush  }\def\:GridPt#1#2#3{   \:xxxx#1
   \advance#1 -#2  \:divide#1#3
   #1\:InCons#1#3  \advance#1 #2
   \ifdim #1=\:xxxx
   \else  \ifdim \:xxxx>#2 \advance#1 #3 \fi \fi  }
\newcount\:IntId    \newcount\:IntCount
\newcount\:DecId    \newcount\:DecCount

\Define\:NewCount{\alloc@ 0\count \countdef \insc@unt }
\Define\:NewDimen{\alloc@ 1\dimen \dimendef \insc@unt }

\def\:NewVar#1#2#3#4#5{ \:multid#1
   \def\:temp{    \csname \string#4\the#4\endcsname\z@
      \edef#1{\noexpand#3  \csname \string#4\the#4\endcsname}}
   \def\:next{  \global#5#4   \expandafter
       #2   \csname \string#4\the#4\endcsname  \:temp }
   \advance#4  \@ne
   \ifnum #4 > #5 \else \def\:next{\:temp}  \fi   \:next}

\def\IntVar#1{\:NewVar#1\:NewCount\:IntOp\:IntId\:IntCount}
\def\DecVar#1{\:NewVar#1\:NewDimen\:DecOp\:DecId\:DecCount}

\DecVar\Q  \DecVar\R   \DecVar\T
\IntVar\I  \IntVar\J   \IntVar\K
\def\WriteVal#1{\immediate\write\sixt@@n{...\string#1=#1;}}

\newdimen\:X   \newdimen\:Y
\newdimen\:x   \newdimen\:y   \newdimen\:d
\newdimen\:xx  \newdimen\:yy  \newdimen\:dd
\newdimen\:xxx \newdimen\:yyy \newdimen\:ddd
\newdimen\:xxxx\newdimen\:yyyy\newdimen\:dddd
\newcount\:J  \newcount\:K
\newdimen\:DI   \newdimen\:DJ
\newdimen\:DK   \newdimen\:DL   \newtoks\:t
\def\:IntFromPt#1#2{
   \:d#2\relax
   \advance\:d  \ifdim\:d<-0.5\p@-\fi  0.5\p@
   #1\:d    \divide#1  65536\relax}
\def\:temp{\catcode`\p12  \catcode`\t12}
\def\:Cons{\catcode`\p11  \catcode`\t11}
\:temp  \def\:Frac#1pt{#1}
        \def\:rnd#1.#2pt{#1}  \:Cons

\def\:Cons#1{\expandafter\:Frac\the#1}
\def\:sqr#1{#1\expandafter\:Frac\the#1#1}
\def\:InCons#1{\expandafter\:rnd\the#1}\def\:Val#1{#1;}
\let\Val\:Val\def\:IntOp#1#2{\csname :Op#2\endcsname#1}
\let\:SvIntOp\:IntOp
\def\:PreIntOp{\let\:IntOp\empty
   \let\Val\empty}
\def\:PostIntOp{\let\:IntOp\:SvIntOp
   \let\Val\:Val}

\expandafter\def\csname :Op;\endcsname#1{ \the#1}
\expandafter\def\csname :Op=\endcsname#1#2;{
   \:PreIntOp#1#2\:PostIntOp}
\expandafter\def\csname :Op+\endcsname#1#2;{
   \:PreIntOp\advance #1  #2\:PostIntOp}
\expandafter\def\csname :Op-\endcsname#1#2;{
   \:PreIntOp\advance #1  -#2\:PostIntOp}
\expandafter\def\csname :Op/\endcsname#1#2;{
   \:PreIntOp\divide#1   #2\:PostIntOp}
\expandafter\def\csname :Op*\endcsname#1#2;{
   \:PreIntOp\multiply#1   #2\:PostIntOp}
\def\:DecOp#1#2{ \csname :xOp#2\endcsname#1}
\let\:SvDecOp\:DecOp
\def\:PreDecOp{\let\:IntOp\the \def\:DecOp{\:Cons}
   \let\Val\empty   \let\:du\empty}
\def\:PostDecOp{\let\:IntOp\:SvIntOp \let\Val\:Val
   \let\:DecOp\:SvDecOp  \let\:du\::du  }

\def\::du#1{\p@
   \ifx#1\p@ \let\:temp\relax
   \else     \def\:temp{\du{#1}}
   \fi\:temp}                    \:PostDecOp

\expandafter\def\csname :xOp;\endcsname#1{ \:Cons#1}
\expandafter\def\csname :Op[\endcsname#1#2];{
   \:PreDecOp \:dd#2\p@  \:IntFromPt#1\:dd
                                 \:PostDecOp  }
\expandafter\def\csname :xOp=\endcsname#1#2;{
   \:PreDecOp#1#2\p@\:PostDecOp               }
\expandafter\def\csname :xOp(\endcsname#1#2){
   \:PreDecOp#1#2\p@\:PostDecOp               }
\expandafter\def\csname :xOp+\endcsname#1#2;{
   \:PreDecOp\advance #1  #2\p@\:PostDecOp  }
\expandafter\def\csname :xOp-\endcsname#1#2;{
   \:PreDecOp\advance #1  -#2\p@\:PostDecOp }
\expandafter\def\csname :xOp*\endcsname#1#2;{
   \:PreDecOp#1 #2#1\:PostDecOp              }
\expandafter\def\csname :xOp/\endcsname#1#2;{
   \:PreDecOp  \:divide#1{#2\p@}  \:PostDecOp }
\let\IF\ifnum

\def\EqText(#1,#2){
   \z@=\z@ \fi  \def\:temp{#1}
                \def\:next{#2}    \ifx \:temp\:next }

\def\:IfInt#1(#2,#3){ \z@=\z@ \fi
   \:IntOp\:K=#2;  \:IntOp\:J=#3; \ifnum  \:K#1\:J }

\def\:IfDim#1(#2,#3){ \z@=\z@ \fi
   \:DecOp\:d=#2;  \:DecOp\:dd=#3; \ifdim  \:d#1\:dd }

 \def\Do(#1,#2)#3{
   \expandafter\let
   \csname :Back\the\:level\endcsname\:Do
\expandafter\edef\csname :DoVars\the\:level\endcsname{
   \DoReg\the\DoReg \:DoB\the\:DoB}
\advance\:level  \@ne
   \DoReg#1  \:DoB#2  \relax
   \ifnum \DoReg<\:DoB
      \def\:Do{\ifnum \DoReg>\:DoB
                  \let\:Do\relax
               \else  #3\advance\DoReg  \@ne \fi
               \:Do}
   \else
      \def\:Do{\ifnum \DoReg<\:DoB
                  \let\:Do\relax
               \else  #3\advance\DoReg \m@ne  \fi
               \:Do}
   \fi  \def\:nextdo{ \:Do \advance\:level  \m@ne
\csname :DoVars\the\:level\endcsname
\def\:temp{\let\:Do}
\expandafter\:temp\csname
   :Back\the\:level\endcsname  } \:nextdo}

\let\:Do\relax

\newcount\DoReg   \let\:DoReg\DoReg       \newcount\:DoB
 \newcount\:level\def\::divide#1{   \:DI\:DK   \:dddd\:DL
   \advance\:DI -\:Cons\:dddd#1
   \:IntFromPt\:J\:dddd  \advance\:dddd -\:J\p@
   \multiply\:J  \@M   \:IntFromPt\:K{\@M\:dddd}
   \advance\:J \:K     \:dddd\@M\p@
   \divide\:dddd \:J   \advance#1 \:Cons\:DI\:dddd  }

\def\:divide#1#2{   \:DK#1   \:DL#2   #1\z@
   \::divide#1  \::divide#1  \::divide#1
   \::divide#1  \::divide#1  }
\def\:Sqrt#1{ \ifdim #1<\:mmp   #1\z@  \else
   \:dd#1   \divide\:dd \tw@
   \def\::Sqrt{  \:ddd#1
      \:divide\:ddd\:dd      \:AbsDif\:d\:dd\:ddd
      \advance\:dd \:ddd   \divide\:dd \tw@
      \ifdim  \:d < \:mmmp
         \let\::Sqrt\relax  \fi
      \::Sqrt}
   \::Sqrt   #1\:dd   \fi }\Define\:length{
   \:dd \:AbsVal \:ddd
   \:abs\:x  \:abs\:y
   \ifdim \:dd<\:x \:dd\:x \fi
   \ifdim \:dd<\:y \:dd\:y \fi
   \ifdim \:dd>\:mp
      \:divide\:x\:dd     \:sqr\:x
      \:divide\:y\:dd     \:sqr\:y
      \:divide\:ddd\:dd   \:sqr\:ddd
      \advance\:x  \:y  \advance\:x \:ddd
      \:y\:dd  \:Sqrt\:x  \:d\:Cons\:x\:y
   \else \:d\:dd \fi}\Define\:distance(2){\MarkLoc(@^)
   \MoveToLoc(#1)  \:x \:X  \:y \:Y  \:ddd\:Z
   \MoveToLoc(#2)
   \advance\:x  -\:X   \advance\:y  -\:Y
   \advance\:ddd  -\if:IIID \:Z \else \:ddd \fi
   \:length  \MoveToLoc(@^)}\def\:NormalizeDeg#1{
   \:DL#1   \:K\:InCons\:DL
   \divide\:K  \:cccvx   \multiply\:K  \:cccvx
   \advance #1 -\:K\p@
   \ifdim #1<\z@ \advance #1  \:CCCVX \fi
   \ifdim #1=\z@
      \ifdim\:DL=\z@ \else
         \advance #1  \:CCCVX \fi \fi }\def\:CosSin#1{ \:DK#1
   \:NormalizeDeg\:DK \def\:tempA{}
\ifdim \:CVXXX<\:DK
   \def\:tempA{\:y-\:y}
   \advance\:DK -\:CCCVX  \:DK-\:DK   \fi
\ifdim \:XC<\:DK
   \edef\:tempA{\:x-\:x \:tempA}
   \advance\:DK -\:CVXXX  \:DK-\:DK   \fi
\ifdim 45\p@<\:DK
   \edef\:tempA{\:d\:x \:x\:y \:y\:d \:tempA}
   \advance\:DK -\:XC   \:DK-\:DK   \fi
   \:x\p@   \:y0.01745\:DK   \:d\:y   \:K\@ne
   \edef\:next{\advance\:K \@ne
      \:sqr\:d  \divide\:d \:K  \advance}
   \:next \:x -\:d   \:next \:y -\:d
   \:next \:x  \:d   \:next \:y  \:d
   \:next \:x -\:d   \:next \:y -\:d
   \:next \:x  \:d   \:next \:y  \:d
   \:tempA   }   \Define\:rInitOval(2){
   \:Zunits\p@   \:dd#1\:Zunits
   \:d\:Cons\:dd\:Xunitsx   \edef\:Xx{\:Cons\:d}
   \:d\:Cons\:dd\:Xunitsy   \edef\:Xy{\:Cons\:d}
   \:dd#2\:Zunits
   \:d\:Cons\:dd\:Yunitsx   \edef\:Yx{\:Cons\:d}
   \:d\:Cons\:dd\:Yunitsy   \edef\:Yy{\:Cons\:d}}
\Define\:xyInitOval(2){
   \:d#1\:Xunits   \edef\:Xx{\:Cons\:d} \def\:Xy{0}
   \:d#2\:Yunits   \edef\:Yy{\:Cons\:d} \def\:Yx{0} }
 \def\:FigSize#1#2#3{
   \:x\:LBorder   \:y\:RBorder   \:d\:TeXLoc
   {\Object\:temp{#3}
    \setbox\:box\hbox{ \:temp
       \multiply\:x by \tw@  \multiply\:y by \tw@
       \xdef\:FSize{ \noexpand#1=\:Cons\:x;
                       \noexpand#2=\:Cons\:y;}}}
   \global\:LBorder\:x   \global\:RBorder\:y
   \global\:TeXLoc \:d
   \:FSize}
\expandafter\let \csname 0:Ln \endcsname\:Ln
\expandafter\def\csname 1:Ln \endcsname{
      \advance\:x -\:X  \advance\:y -\:Y
      \csname 0:Ln \endcsname(\:x,\:y)  }

\newcount\:ClipLevel   \:ClipLevel\@ne

\Define\Clip{\futurelet\:next\:Clip}

\Define\:Clip{
   \ifx \:next[  \expandafter\:DefClipOut
        \else    \expandafter\:DefClip     \fi }

\def\:DefClipOut[#1]{ \:DefClip(#1) }
\Define\:DefClip(1){  \def\:temp{#1}
   \ifx \:temp\empty
      \:ClipLevel\@ne     \def\:next{\let\:Ln}
      \expandafter\:next\csname 0:Ln \endcsname
   \else  \def\:temp{\::DefClip(#1)}  \fi  \:temp }
\Define\::DefClip(1){  \MarkLoc(^)
   \:x\:X  \:y \:Y  \Move(#1)
   \ifdim\:x>\:X  \:dd\:X  \:X\:x  \:x\:dd  \fi
   \ifdim\:y>\:Y  \:dd\:Y  \:Y\:y  \:y\:dd  \fi
   \advance\:ClipLevel  \@ne
   \expandafter\edef\csname \the
      \:ClipLevel :Ln \endcsname{
          \:xxx \the\:x   \:yyy \the\:y
          \:xxxx\the\:X  \:yyyy\the\:Y
          \ifx \:next[   \noexpand\:ClipOut
          \else          \noexpand\:ClipIn     \fi }
   \let\:Ln\:ClipLn   \MoveToLoc(^)  }\def\:ClipLn(#1,#2){
   \:x#1\:Xunits \:y#2\:Yunits
   {  \let\:Xunits\empty  \let\:Yunits\empty
   \advance\:x \:X   \advance\:y \:Y
   \ifdim \:x<\:X  \:dd\:X \:X\:x \:x\:dd
                   \:dd\:Y \:Y\:y \:y\:dd  \fi
   \:diff\:dd\:X\:x   \:diff\:ddd\:Y\:y
\:Z \:AbsVal \:dd
\advance\:Z  \:AbsVal\:ddd
\ifdim \:Z>\sixt@@n\p@
   \divide\:dd   128
   \divide\:ddd  128  \fi
\:Z\:Cons\:y\:dd
\advance\:Z -\:Cons\:x\:ddd
\ifdim \:dd<\z@
  \:dd-\:dd  \:ddd-\:ddd  \:Z-\:Z
\fi
   \csname \the \:ClipLevel :Ln \endcsname      }
   \advance\:X \:x   \advance\:Y \:y  }\Define\:ClipIn{
   \def\:next{\let\:next}
   \expandafter\:next\csname \the
      \:ClipLevel :Ln \endcsname
   \advance\:ClipLevel \m@ne
   { \:ClipLeft\:xxxx \:ClipDown\:yyy \:ClipUp\:yyyy \:next }
   { \:ClipRight\:xxx \:ClipDown\:yyy \:ClipUp\:yyyy \:next  }
   { \:ClipDown\:yyyy \:next }
     \:ClipUp\:yyy    \:next }

\Define\:ClipOut{
   \:ClipLeft\:xxx      \:ClipRight\:xxxx
   \:ClipUp  \:yyyy     \:ClipDown \:yyy
   \def\:next{\let\:next}
   \expandafter\:next
      \csname \the\:ClipLevel :Ln \endcsname
   \advance\:ClipLevel \m@ne      \:next  }\def\:ClipLeft#1{
   \ifdim       \:x<#1  \:KilledLine
   \else \ifdim \:X<#1  \:X#1
      \ifdim \:dd>\:mmmp
         \:Y\:Cons\:ddd\:X  \advance\:Y \:Z
         \:divide\:Y\:dd
   \fi \fi \fi     \:CondKilLn  }

\def\:ClipRight#1{
   \ifdim       \:X>#1  \:KilledLine
   \else \ifdim \:x>#1  \:x#1
      \ifdim \:dd>\:mmmp
         \:y\:Cons\:ddd\:x  \advance\:y \:Z
         \:divide\:y\:dd
   \fi \fi \fi    \:CondKilLn }

\Define\:CondKilLn{
   \:d\:x  \advance\:d -\:X
   \ifdim \:d<\z@  \:d-\:d \fi
   \ifdim \:y<\:Y  \advance\:d \:Y \advance\:d -\:y
   \else           \advance\:d \:y \advance\:d -\:Y  \fi
   \ifdim  \:d<\:mmp   \:KilledLine \fi
   \ifdim  \:thickness=\z@ \:KilledLine \fi}

\Define\:KilledLine{
   \let\:ClipLeft\:gobble  \let\:ClipRight\:gobble
   \let\:ClipUp  \:gobble  \let\:ClipDown \:gobble
   \let\:next\relax
   \expandafter\def\csname 1:Ln \endcsname{}}\def\:ClipUp#1{
   \:AbsDif\:d\:y\:Y
   \ifdim  \:d<\:ragged
      \advance\:y  0.5\:thickness
\advance\:Y -0.5\:thickness
\ifdim       \:Y>#1  \:KilledLine
\else \ifdim \:y>#1
   \:thickness#1 \advance\:thickness -\:Y
   \advance\:Y  0.5\:thickness  \:y\:Y
\else
   \advance\:Y  0.5\:thickness
   \advance\:y -0.5\:thickness
\fi  \fi
\:dd\p@   \:ddd\z@
\def\:temp{  \:Z\:Y }  \:temp
   \else \let\:temp\relax
          \ifdim \ifdim\:Y<\:y\:Y\else\:y\fi >#1  \:KilledLine
   \else  \ifdim  \::ClipUp#1\:X\:Y
   \else  \ifdim  \::ClipUp#1\:x\:y
   \fi \fi \fi \fi   \:CondKilLn  }

\def\::ClipUp#1#2#3{
#3>#1   #3#1
\ifdim \:AbsVal\:ddd>\:mmmp
   #2\:Cons\:dd#3  \advance#2 -\:Z
   \:divide#2\:ddd
\fi  }\def\:ClipDown#1{   \:Ex2#1
   \:Flip\:y  \:Flip\:Y  \:ClipUp#1
   \:Flip\:y  \:Flip\:Y  \:temp }

\def\:Flip#1{   #1-#1  \advance#1 \:Ex  }
\Define\:SetDrawWidth{
   \hsize\:RBorder      \advance\hsize -\:LBorder
   \leftskip -\:LBorder \rightskip\z@}
\newdimen\:thickness   \:thickness0.75\p@

\Define\PenSize(1){\:thickness#1\relax}
\Define\:Draw{   \ifvmode \noindent\hfil\fi
   \global\:TeXLoc\z@
   \vbox\bgroup
           \begingroup
   \def\EndDraw{
           \endgroup   \:SetDrawWidth
         \egroup}
   \:DraCatCodes      \parindent\z@    \everypar{}
\leftskip\z@     \rightskip\z@    \boxmaxdepth\maxdimen
\let\FigSize\:FigSize
\def\Draw{\:wrn1{}} \:CommonIID   \:InDraw }
\Define\:CommonIID{\def\RotateTo{\:RotateTo}

\def\MoveF{\:MvF}\def\LineToLoc{\:LnToLoc}}\newdimen\:Xunits   \:Xunits\p@
\newdimen\:Yunits   \:Yunits\p@
\Define\:InDraw{\:Opt()\::InDraw{\:Xunits,\:Yunits}}
\Define\::InDraw(2){   \:X\z@   \:Y\z@
   \:Xunits#1 \relax  \:Yunits#2  \:AdjRunits
   \global\:LBorder\z@    \global\:RBorder\z@
   \:loadIID   \leavevmode }

\Define\:Resize(2){
   \:Xunits#1\:Xunits  \relax
   \:Yunits#2\:Yunits  \:AdjRunits  }
\Define\:Units(2){
   \:Xunits#1 \relax   \:Yunits#2    \:AdjRunits   }
\Define\:loadIID{
   \def\LineAt{\:DLn}
   \def\LineTo{\:LnTo}
   \def\MoveTo{\:MvTo}
   \def\Line{\:Ln}
   \def\Move{\:Mv}
   \def\MoveF{\:MvF}
   \:rotatedfalse\def\:InitOval{\:xyInitOval}
   
   \def\Units{\:Units}}

\def\DrawOn{\def\Draw{\:Draw}}                      \DrawOn
\def\DrawOff{\def\Draw{\begingroup \:J\@cclv
                       \:NoDrawSpecials \:NoDraw}}
\catcode`\/0 \catcode`\\11
/def/:NoDraw#1\EndDraw{/endgroup}
/catcode`/\0 /catcode`//12

\def\:NoDrawSpecials{\catcode\:J11
  \ifnum \:J=\z@
     \let \:NoDrawSpecials\relax \fi
  \advance\:J  \m@ne \:NoDrawSpecials}
\let\:XunitsReg\:Xunits   \let\:YunitsReg\:Yunits  \Define\EntryExit(4){
   \edef\:InOut##1{
      \noexpand\ifcase ##1\space
         #1\noexpand\or #2\noexpand\or
         #3\noexpand\or #4\noexpand\fi}}

\EntryExit(0,0,0,0)\Define\:DrawBox{   \:x0.5\wd\:box
   \:y\ht\:box \advance\:y  \dp\:box
   \divide\:y \tw@
   \advance \:X -\:InOut0\:x
   \advance \:Y -\:InOut1\:y
   { \advance \:X -\:x   \advance \:Y -\:y
     { \ifdim\:X<\:LBorder \global\:LBorder\:X\fi
\advance\:X \wd\:box
\ifdim\:X>\:RBorder \global\:RBorder\:X\fi }
     \advance \:Y  \dp\:box
     {\:d\:X \advance\:d \wd\:box
 \advance\:X -\:TeXLoc   \global\:TeXLoc\:d
 \vrule width\:X depth\z@ height\z@
 \raise \:Y \box\:box} }
   \edef\MoveToExit(##1,##2){
   \:X\the\:X   \:Y\the\:Y
   \:x\the\:x   \:y\the\:y
   \advance\:X  ##1\:x
   \advance\:Y  ##2\:y}
   \advance \:X  \:InOut2\:x
   \advance \:Y  \:InOut3\:y }
 \Define\ThreeDim{\:Opt[]\:ThreeDim{\p@}}

\def\:ThreeDim[#1](#2){\::ThreeDim[#1](#2,,)}

\def\::ThreeDim[#1](#2,#3,#4,#5){ \bgroup\begingroup
   \def\EndThreeDim{          \endgroup\egroup}
   \:IIIDtrue   \:Zunits#1  \:Z\z@
   \def\:temp{#4}
   \ifx \:temp\empty   \:CosSin{#3\p@}  \:divide\:x\:y       \:Ey\:x
\:CosSin{#2\p@}  \:Ex\:Cons\:Ey\:x  \:Ey\:Cons\:Ey\:y
\let\:project\:projectPar

   \else               \:Ex#2\:Xunits \:Ey#3\:Yunits \:Ez#4\:Zunits
\let\:project\:projectPer
 \fi
   \def\LineAt{\:tDLn}
\def\LineTo{\:tLnTo}
\def\MoveTo{\:tMvTo}
\def\Line{\:tLn}
\def\Move{\:tMv}

\def\Units{\:tUnits}\def\RotateTo{\:tRotateTo}

\def\MoveF{\:tMvF}
\:Vdirection\z@ \def\LineToLoc{\:tLnToLoc} }  \Define\:projectPer{
   \:diff\:x\:X\:Ex   \:diff\:y\:Y\:Ey
   \:diff\:xxxx\:Z\:Ez
   \:divide\:x\:xxxx      \:divide\:y\:xxxx
   \:x-\:Cons\:Ez\:x     \:y-\:Cons\:Ez\:y
   \advance\:x  \:Ex    \advance\:y \:Ey}
\Define\:projectPar{
   \:x\:Cons\:Ex\:Z  \advance\:x \:X
   \:y\:Cons\:Ey\:Z  \advance\:y \:Y}
\def\:tDLn(#1,#2,#3,{\:tMvTo(#1,#2,#3)
                     \:tLnTo(}
\Define\:tMvTo(3){   \:X#1\:Xunits
   \:Y#2\:Yunits    \:Z#3\:Zunits}
\Define\:tMv(3){     \advance\:X  #1\:Xunits
   \advance\:Y  #2\:Yunits
   \advance\:Z  #3\:Zunits}
\Define\:tResize(3){   \:Xunits#1\:Xunits
   \:Yunits#2\:Yunits \:Zunits#3\:Zunits
   \:AdjRunits}
\Define\:tUnits(3){    \:Xunits#1  \relax
   \:Yunits#2  \relax       \:Zunits#3
   \:AdjRunits}\Define\:tLnTo(3){
   \:project   \edef\:temp{\:x\the\:x \:y\the\:y}
   \:X#1\:Xunits  \:Y#2\:Yunits  \:Z#3\:Zunits
   \:DLN}

\Define\:tLn(3){
   \:project   \edef\:temp{\:x\the\:x \:y\the\:y}
   \advance\:X  #1\:Xunits
   \advance\:Y  #2\:Yunits
   \advance\:Z  #3\:Zunits   \:DLN}

\Define\:DLN{
   \TwoDim      \:temp        \LineTo(\:x\du,\:y\du)
   \EndTwoDim }\Define\TwoDim{\bgroup\begingroup
   \def\EndTwoDim{\endgroup\egroup}
   \:loadIID
   \if:IIID  \:IIIDfalse \:project \:X\:x \:Y\:y
             \:CommonIID \fi
   \Units(\:Xunits,\:Yunits)}
\newif\if:rotated
\Define\RotatedAxes(2){\begingroup
   \def\EndRotatedAxes{\endgroup}
   \if:IIID  \:IIIDfalse
             \:project \:X\:x \:Y\:y
             \:CommonIID    \fi
   \:DK#2\p@ \advance\:DK -#1\p@
\advance\:DK -\:CVXXX   \:NormalizeDeg\:DK
\ifdim \:DK<\:mmp \:wrn4{(#1,#2)}
\else
  \:DK#1\p@ \:NormalizeDeg\:DK    \:CosSin\:DK
  \:Xunitsx\:x  \:Xunitsy\:y
  \:DK#2\p@ \advance\:DK -\:XC \:NormalizeDeg\:DK \:CosSin\:DK
  \edef\Units(##1,##2){\noexpand\:Units(##1,##2)
    \:Xunitsx \:Cons\:Xunitsx\:Xunits
    \:Xunitsy \:Cons\:Xunitsy\:Xunits
    \:Yunitsx-\:Cons\:y\:Yunits
    \:Yunitsy \:Cons\:x\:Yunits  } \fi
 \def\MoveTo{\:rMvTo}
\def\Move{\:rMv}
\def\LineTo{\:rLnTo}
\def\Line{\:rLn}
\def\MoveF{\:rMvF}\def\:InitOval{\:rInitOval}

   \MarkLoc(:org) \Units(\:Xunits,\:Yunits)  \:rotatedtrue }
\newdimen\:Xunitsx  \newdimen\:Xunitsy
\newdimen\:Yunitsx  \newdimen\:Yunitsy
\Define\:rResize(2){
   \:Xunits #1\:Xunits    \:Yunits #2\:Yunits
   \:Xunitsx#1\:Xunitsx   \:Xunitsy#1\:Xunitsy
   \:Yunitsx#2\:Yunitsx   \:Yunitsy#2\:Yunitsy  }
\Define\:rMvTo{\MoveToLoc(:org) \:rMv}
\Define\:rMv(1){ \:rxy(#1)
   \advance\:X  \:x \advance\:Y  \:y}
\Define\:rLnTo(1){ \MarkLoc(:a) \:rMvTo(#1) {\LineToLoc(:a)} }
\Define\:rLn(1){ \:rxy(#1) \:Ln(\:x\du,\:y\du) }
\Define\:rxy(2){ \:Zunits\p@   \:Zunits#1\:Zunits
   \:x\:Cons\:Zunits\:Xunitsx
   \:y\:Cons\:Zunits\:Xunitsy
   \:Zunits\p@   \:Zunits#2\:Zunits
   \advance\:x  \:Cons\:Zunits\:Yunitsx
   \advance\:y  \:Cons\:Zunits\:Yunitsy}
\Define\:rMvF(1){  \:CosSin\:direction
   \:Zunits\p@            \:Zunits#1\:Zunits
   \:x\:Cons\:Zunits\:x   \:y\:Cons\:Zunits\:y
   \edef\:temp{(\:Cons\:x,\:Cons\:y)}   \expandafter\Move\:temp}
\Define\:AdjRunits{
   \:Xunitsx\:Xunits   \:Xunitsy\z@
   \:Yunitsx\z@        \:Yunitsy\:Yunits}
\Define\Text{  \setbox\:box
   \vtop\bgroup    \edef\DoReg{\the\DoReg}
      \hyphenpenalty\@M  \exhyphenpenalty\@M
      \catcode`\ 10 \catcode`\^^M13 \catcode`\^^I10
      \catcode`\&4  \let~\space
      \:Text}                       \catcode`\^^M13 %
\def\:Text(--#1--){%
      \:SetLines#1\hbox{}^^M--)^^M %
   \egroup                  %
   \if:IIID \TwoDim  \:DrawBox  \EndTwoDim %
   \else             \:DrawBox  \fi}       %
\def\:SetLines#1^^M{        %
   \def\:TextLine{#1}       %
   \ifx \:TextLine\:LastLine   \let\:temp\relax      %
   \else  \def\:temp{                                 %
             \:IndirectLines#1\relax~~--)~~\:SetLines}%
   \fi  \:temp }                      \catcode`\^^M9
\def\:IndirectLines#1~~{    \def\:TextLine{#1}
  \ifx \:TextLine\:LastLine   \let\:temp\relax
  \else  \def\:temp{\:AddLine{#1}\:IndirectLines}
  \fi  \:temp }

\Define\:LastLine{--)}

\def\:AddLine#1{
   \ifvmode \noindent    \hsize\z@ \else
      \hfil \penalty-500 \hbox{}    \fi
   \hfil#1
   \setbox\:box\hbox{#1}
   \ifdim \wd\:box>\hsize \hsize\wd\:box \fi}\def\TextPar#1#2{
   \def\:TxtPar##1(##2){##1(--##2--)}
   \edef\:temp{\expandafter\noexpand\csname :\string#2\endcsname}
   \edef#2{\noexpand\:TextPar\expandafter\noexpand\:temp}
   \expandafter\let\:temp\:undefined
   \expandafter#1\:temp}
                                   \catcode`\^^M13
\def\:TextPar#1{\begingroup     \catcode`\&4         %
   \catcode`\ 10 \catcode`\^^M13 \catcode`\^^I10 %
   \:TPar{#1}}                     \catcode`\^^M9  %

\def\:TPar#1(--#2--){\endgroup
   #1(--#2--)  }
 \newdimen\:direction  \newdimen\:Vdirection
\Define\:Rotate(1){   \advance \:direction   #1\p@
                      \:NormalizeDeg\:direction  }
\Define\:RotateTo(1){ \:direction  #1\p@
                      \:NormalizeDeg\:direction  }

\Define\:tRotate(2){
   \advance\:Vdirection  #2\p@
   \:NormalizeDeg\:Vdirection  \:Rotate(#1)}
\Define\:tRotateTo(2){
   \:Vdirection #2\p@
   \:NormalizeDeg\:Vdirection  \:RotateTo(#1)}
\Define\LineF(1){\MarkLoc(,)\MoveF(#1) {\LineToLoc(,)}}
\Define\:MvF(1){    \:CosSin\:direction
   \:d#1\:Xunits   \advance\:X  \:Cons\:d\:x
   \:d#1\:Yunits   \advance\:Y  \:Cons\:d\:y   }
\Define\MoveFToOval(2){
   \:NormalizeDeg\:direction \:CosSin\:direction  \:Zunits\p@
   \:d#2\:Zunits    \:x\:Cons\:d\:x
   \ifdim \:d=\z@  \:err5{...,\:Cons\:d} \fi
   \:d#1\:Zunits    \:y\:Cons\:d\:y
   \ifdim \:d=\z@  \:err5{\:Cons\:d,...} \fi     \:SearchDir
   \XSaveUnits
      \if:rotated
        \:Zunits\p@     \:Zunits#1\:Zunits
        \:x\:Cons\:Zunits\:Xunits
        \:Zunits\p@     \:Zunits#2\:Zunits
        \:y\:Cons\:Zunits\:Yunits
      \else
        \:x#1\:Xunits  \:y#2\:Yunits  \fi     \relax
      \Units(\:x,\:y)
      \:Vdirection\:direction  \:direction\:ddd
      \MoveF(\@ne)  \:direction\:Vdirection  \XRecallUnits  }
\Define\:tMvF(1){    \:CosSin\:Vdirection
   \:d #1\:Zunits   \advance\:Z  \:Cons\:y\:d
   \:xx#1\:Xunits     \:yy#1\:Yunits
   \:xx\:Cons\:x\:xx    \:yy\:Cons\:x\:yy
   \:CosSin\:direction
   \:x\:Cons\:xx\:x    \:y\:Cons\:yy\:y
   \advance\:X  \:x   \advance\:Y  \:y } \Define\CSeg{\:Opt[]\:CSeg1}
\def\:CSeg[#1]#2(#3,#4){   \MarkLoc($^)
   \MoveToLoc(#4) \:x\:X \:y\:Y
   \if:IIID  \:d\:Z  \fi   \MoveToLoc(#3)
   \advance\:x -\:X  \:x#1\:x
   \advance\:y -\:Y  \:y#1\:y
   \if:IIID   \advance\:d -\:Z  \:d#1\:d \fi
   \:t{#2}
   \edef\:temp{\the\:t(
      \expandafter\:Frac\the\:x\noexpand\:du,
      \expandafter\:Frac\the\:y\noexpand\:du \if:IIID ,
      \expandafter\:Frac\the\:d\noexpand\:du \fi)}
   \MoveToLoc($^)     \:temp}\Define\LSeg{\:Opt[]\:LSeg1}
\def\:LSeg[#1]#2(#3,#4){   \:distance(#3,#4)
   \:d#1\:d  \:t{#2}
   \edef\:temp{\the\:t(\expandafter\:Frac\the\:d\noexpand\:du)}
   \:temp}\Define\DSeg{\:Opt[]\:DSeg1}

\def\:DSeg[#1]#2(#3,#4){   \MarkLoc(^)
   \MoveToLoc(#4)  \:xxx\:X   \:yyy\:Y \:xxxx\:Z
   \MoveToLoc(#3)
   \advance\:xxx -\:X   \advance\:yyy -\:Y
   \ifdim \:AbsVal\:xxx<\:mmmp
      \ifdim \:AbsVal\:yyy<\:mmmp \:wrn5{#3,#4}
   \fi\fi
   \if:IIID
      \advance \:xxxx  -\:Z
      \:divide\:xxx\:Xunitsx
      \:divide\:yyy\:Yunitsy
      \:divide\:xxxx\:Zunits
      \:x\:xxx  \:y\:yyy   \:ddd\z@  \:length
      \:x\:d    \:y\:xxxx  \:SearchDir
      \:yyyy\:ddd
      \:x\:xxx  \:y\:yyy
   \else
      \:x  \:AbsVal\:Yunitsy
\:y  \:AbsVal\:Xunitsy  \ifdim \:y>\:x  \:x\:y  \fi
\:y  \:AbsVal\:Yunitsx  \ifdim \:y>\:x  \:x\:y  \fi
\:y  \:AbsVal\:Xunitsx  \ifdim \:y>\:x  \:x\:y  \fi
\:K  \:InCons\:x  \relax
      \ifnum \:K<\thr@@     \:K\@ne
\else \ifnum \:K<\sixt@@n   \:K4
\else \ifnum \:K<\:XC       \:K\sixt@@n
\else \ifnum \:K<\@m        \:K\@cclvi
\fi \fi \fi \fi
\divide\:xxx \:K
\divide\:yyy \:K
      \:x \:Cons\:Yunitsy\:xxx
\advance\:x -\:Cons\:Yunitsx\:yyy
      \:y-\:Cons\:Xunitsy\:xxx
\advance\:y \:Cons\:Xunitsx\:yyy
   \fi
   \:SearchDir   \:ddd#1\:ddd   \:t{#2}
   \edef\:temp{\the\:t(
      \:Cons\:ddd \if:IIID ,\:Cons\:yyyy \fi)}
   \MoveToLoc(^) \:temp}   \def\:theDoReg{\def\DoReg{\the\:DoReg}}

\Define\MarkLoc{  \:theDoReg
   \expandafter\edef \csname \:MarkLoc}

\Define\MarkGLoc{  \:theDoReg
   \expandafter\xdef \csname \:MarkLoc}

\Define\:MarkLoc(1){ Loc\space#1:\endcsname{
   \:X\the\:X  \:Y\the\:Y
   \if:IIID \:Z\the\:Z \fi  }       \let\DoReg\:DoReg }

\Define\MoveToLoc(1){   \:theDoReg   \expandafter\ifx
   \csname Loc\space#1:\endcsname\relax \:err2{#1}\fi
   \csname Loc\space#1:\endcsname     \let\DoReg\:DoReg  }

\Define\MarkPLoc(1){  \:theDoReg
   \if:IIID   \:project
      \expandafter\edef \csname Loc\space#1:\endcsname{
         \:X\the\:x  \:Y\the\:y  \:Z\z@}
   \else \:err1\MarkPLoc  \fi   \let\DoReg\:DoReg}

\Define\WriteLoc(1){{ \:theDoReg \edef\:temp{#1}
   \ifx \:temp\empty \else \MoveToLoc(#1) \fi
   \immediate\write\sixt@@n{...
      \:temp=(\the\:X,\the\:Y\if:IIID,\the\:Z\fi)}}}

\Define\:LnToLoc(1){
   \:x\:X \:y\:Y   \MoveToLoc(#1)
   { \:LnTo(\:x\du,\:y\du) }}

\Define\:tLnToLoc(1){
   \:xx\:X \:yy\:Y \:dd\:Z     \MoveToLoc(#1)
   { \:tLnTo(\:xx\du,\:yy\du,\:dd\du) }}\def\:GetLine#1#2#3#4#5{
   \MoveToLoc(#1)   \divide\:X \:eight  \divide\:Y \:eight
   #3\:X  #4\:Y
   \MoveToLoc(#2)   \divide\:X \:eight  \divide\:Y \:eight
   \advance #3 -\:X  \advance #4 -\:Y
   #5\:Cons#3\:Y      \advance #5 -\:Cons#4\:X
   \divide #3 \:eight     \divide  #4 \:eight \relax      }
\def\MoveToLL(#1,#2)(#3,#4){
   \:GetLine{#1}{#2}\:x \:y \:xxx
   \:GetLine{#3}{#4}\:xx\:yy\:xxxx
   \:ddd \:Cons\:x \:yy     \advance\:ddd -\:Cons\:xx\:y
   \ifdim  \:AbsVal\:ddd < \:mmmp
      \:X\@cclv\p@  \:Y\:X
      \:wrn3{(\string#1,\string#2)(\string#3,\string#4)}
   \else
      \:divide\:xxx\:ddd    \:divide\:xxxx\:ddd
      \:X\:Cons\:xxx\:xx   \advance\:X -\:Cons\:xxxx\:x
      \:Y\:Cons\:xxx\:yy   \advance\:Y -\:Cons\:xxxx\:y
   \fi   }\Define\MoveToCC{\:Opt[]\:MoveToCC{}}

\def\:MoveToCC[#1](#2,#3)(#4,#5){
  \:UserUnits(#2,#3)(#4,#5)
  \:distance($#2,$#4)
\ifnum \:d<\:mp \MoveToLoc(#3)
   \:wrn3{(#1,#2)(#3,#4)}
\else                 \:xx \:d
   \:distance($#2,$#3)  \:xxx\:d
   \:distance($#5,$#4)
     \:yy \:xxx  \advance\:yy -\:d
\:yyy\:xxx  \advance\:yyy \:d
\:divide\:yy\:xx     \:yy\:Cons\:yy\:yyy
\advance\:yy \:xx  \divide\:yy \tw@
     \:yyy \ifdim \:AbsVal\:xxx>\:AbsVal\:yy \:xxx \else \:yy \fi
\ifdim \:AbsVal\:yyy<\:mp \:yyy\z@ \else
   \:divide\:xxx\:yyy  \:sqr\:xxx
   \:yyyy\:yy  \:divide\:yyyy\:yyy  \:sqr\:yyyy
   \advance\:xxx -\:yyyy   \:Sqrt\:xxx
   \:yyy\:Cons\:yyy\:xxx
\fi
     \MoveToLoc($#4)  \:x\:X  \:y\:Y  \:divide\:yy\:xx
\MoveToLoc($#2)
\advance\:x -\:X     \advance\:y -\:Y
\advance\:X \:Cons\:yy\:x
\advance\:Y \:Cons\:yy\:y
     \:divide\:yyy\:xx
\advance\:X  #1\:Cons\:yyy\:y
\advance\:Y -#1\:Cons\:yyy\:x

  \fi  \:SysUnits  }\def\:UserUnits(#1,#2)(#3,#4){
   \:xx\:Xunitsx  \:xxx\:Xunitsy
   \:yy\:Yunitsx  \:yyy\:Yunitsy
   \:xxxx\:Cons\:yy\:xxx  \advance\:xxxx \:Cons\:yyy\:xx
   \ifdim \:AbsVal\:xxxx>\:mmmp
      \:divide\:xx\:xxxx   \:divide\:xxx{-\:xxxx}
      \:divide\:yy{\:xxxx} \:divide\:yyy\:xxxx
   \fi
   \:UnLoc(#1)  \:UnLoc(#2)
   \:UnLoc(#3)  \:UnLoc(#4)}\Define\:SysUnits{
   \MoveTo(\:Cons\:X,\:Cons\:Y)}\Define\:UnLoc(1){
   \MoveToLoc(#1)  \:d\:Cons\:yyy\:X
   \advance\:d \:Cons\:yy\:Y
   \:Y\:Cons\:xx\:Y
   \advance\:Y  \:Cons\:xxx\:X   \:X\:d
   \MarkLoc($#1)}

\def\:MoveToLC[#1](#2,#3)(#4,#5){
   \:UserUnits(#2,#3)(#4,#5)
   \MoveToLoc($#2)  \:x\:X  \:y\:Y
\MoveToLoc($#3)  \advance\:x -\:X
                \advance\:y -\:Y
   \edef\:temp{ \:xxx\the\:x  \:yyy\the\:y }
\MoveToLoc($#4)
\advance\:X \:y  \advance\:Y -\:x
\MarkLoc(^$) \MoveToLL($#4,^$)($#2,$#3)
\MarkLoc(^$)
   \:distance($#4,$#5)  \:xx\:d  \:distance($#4,^$)
\ifdim      \:d>\:xx        \:wrn3{(#2,#3)(#4,#5)}
\else \ifdim \:d<\:mmp     \:yy\:xx   \else   \:yy\:d
      \:divide\:yy\:xx  \:sqr\:yy
      \:yy-\:yy   \advance\:yy \p@
      \:Sqrt\:yy  \:yy\:Cons\:xx\:yy
\fi \fi
   \:temp   \:x\:xxx  \:y\:yyy  \:length  \:xx\:d
\:divide\:xxx\:xx
\:divide\:yyy\:xx
\advance\:X  #1\:Cons\:yy\:xxx
\advance\:Y  #1\:Cons\:yy\:yyy
   \:SysUnits }  \def\Object#1{\:Opt(){\:DefineSD#1}0}

\def\:DefineSD#1(#2){\begingroup  \:multid#1
   \:DraCatCodes   \:DefSD#1(#2)}

\def\:DefSD#1(#2)#3{
   \expandafter\::Define\csname\string#1.\endcsname(#2){
      \:t{\:SubD{#3}}
      \if:IIID \edef\:temp{\noexpand\TwoDim \the\:t
                           \noexpand\EndTwoDim}
      \else    \def\:temp{\the\:t}   \fi         \:temp}
   \def#1{\def\:SDname{\csname\string#1.\endcsname}
          \:Opt[]\:CallSD{}}}

\def\:CallSD[#1]{ \edef\:Entry{#1} \:SDname }\def\:SubD#1{
   \let\::RecallXLoc\:AddXLoc   \gdef\:AddXLoc{}
   \edef\:RecallBor{ \global\:LBorder  \the\:LBorder
                  \global\:RBorder  \the\:RBorder
                  \global\:TeXLoc\the\:TeXLoc }
\global\:TeXLoc\z@
\setbox\:box\vbox{\EntryExit(0,0,0,0)
\begingroup
   
   \:InDraw  #1
\endgroup
\:SetDrawWidth                    \let\:XLoc\relax
\xdef\:AddXLoc{\:dd\the\:LBorder  \:AddXLoc}}
\:RecallBor
   \:ddd\dp\:box
   \ifx \:Entry\empty
   \:DrawBox
\else
   \let\:RecallIn\:InOut
   \:x\:X    \:y\:Y
   \def\:XLoc(##1,##2,##3){
      \def\:temp{##1}
      \ifx \:temp\:Entry \:X\:x  \advance\:X -##2
                         \:Y\:y  \advance\:Y -##3
      \fi}
   \:AddXLoc
   \advance\:X  \:dd      \advance\:Y -\:ddd
   \EntryExit(-1,-1,\:InOut2,\:InOut3)  \:DrawBox
   \let\:InOut\:RecallIn
\fi
   \MarkLoc(^)
      \MoveToExit(-1,-1)
   \:xxx\:X   \:yyy\:Y   \advance\:yyy  \:ddd
   \def\:XLoc(##1,##2,##3){
      \:X\:xxx  \advance\:X  ##2     \advance\:X -\:dd
      \:Y\:yyy  \advance\:Y  ##3
      \MarkLoc(##1)}
   \:AddXLoc
   \MoveToLoc(^)
   \ifx \:Entry\empty \else     \MoveToLoc(\:Entry) \fi
   \global\let\:AddXLoc\::RecallXLoc }
\Define\:MarkXLoc(1){  \:theDoReg
   \let\:XLoc\relax
   \xdef\:AddXLoc{\:AddXLoc \:XLoc(#1,\the\:X,\the\:Y)}
   \let\DoReg\:DoReg}
   \catcode`\ 10 \catcode`\^^M5 \catcode`\^^I10
   \let\wlog\:wlog  \let\:wlog\:undefined
\:RestoreCatcodes     \tracingstats1

\input epsf
\title{Noncommutative localization in topology}
\author{Andrew Ranicki}

\parskip\medskipamount
\def\Z{\mathbb Z}
\def\R{\mathbb R}
\def\bysame{\leavevmode\hbox to3em{\hrulefill}\space}
\begin{document}
\maketitle

\section*{Introduction}

The topological applications of the Cohn noncommutative localization
considered in this paper deal with spaces (especially manifolds) with
infinite fundamental group, and involve localizations of infinite group
rings and related triangular matrix rings.  Algebraists have usually
considered noncommutative localization of rather better behaved rings,
so the topological applications require new algebraic techniques.

Part 1 is a brief survey of the applications of noncommutative
localization to topology: finitely dominated spaces, codimension 1 and
2 embeddings (knots and links), homology surgery theory, open book
decompositions and circle-valued Morse theory.  These applications
involve chain complexes and the algebraic $K$- and $L$-theory of the
noncommutative localization of group rings.

Part 2 is a report on work on chain complexes over generalized
free products and the related algebraic $K$- and $L$-theory, from
the point of view of noncommutative localization of triangular
matrix rings. Following Bergman and Schofield, a generalized free
product of rings can be constructed as a noncommutative
localization of a triangular matrix ring.  The novelty here is the
explicit connection to the algebraic topology of manifolds with a
generalized free product structure realized by a codimension 1
submanifold, leading to noncommutative localization proofs of the
results of Waldhausen and Cappell on the algebraic $K$- and
$L$-theory of generalized free products.  In a sense, this is more
in the nature of an application of topology to noncommutative
localization! But this algebra has in turn topological
applications, since in dimensions $\geqslant 5$ the surgery
classification of manifolds within a homotopy type reduces to
algebra.

\setcounter{section}{1}
\setcounter{subsection}{0}

\section*{Part 1. A survey of applications}\label{survey}

We start by recalling the universal noncommutative localization of
P.M.Cohn \cite{Co}. Let $A$ be a ring, and let $\Sigma=\{s:P \to Q\}$
be a set of morphism of f.g. projective $A$-modules. A ring morphism
$A \to R$ is {\it $\Sigma$-inverting} if for every $s \in \Sigma$
the induced morphism of f.g. projective $R$-modules
$1 \otimes s : R \otimes_A P \to R \otimes_A Q$ is an isomorphism.
The noncommutative localization $A \to \Sigma^{-1}A$ is
$\Sigma$-inverting, and has the universal property that any
$\Sigma$-inverting ring morphism $A \to R$ has a unique factorization
$A \to \Sigma^{-1}A \to R$. The applications to topology involve
homology with coefficients in a noncommutative localization $\Sigma^{-1}A$.

Homology with coefficients is defined as follows. Let $X$ be a
connected topological space with universal cover $\widetilde{X}$, and let the
fundamental group $\pi_1(X)$ act on the left of $\widetilde{X}$, so
that the (singular) chain complex $S(\widetilde{X})$ is a free left
$\Z[\pi_1(X)]$-module complex.  Given a morphism of rings ${\mathcal
F}:\Z[\pi_1(X)] \to \Lambda$ define the {\it $\Lambda$-coefficient
homology of $X$} to be
$$H_*(X;\Lambda)~=~H_*(\Lambda\otimes_{\Z[\pi_1(X)]}S(\widetilde{X}))~.$$
If $X$ is a $CW$ complex then $S(\widetilde{X})$ is chain equivalent
to the cellular free $\Z[\pi_1(X)]$-module chain complex $C(\widetilde{X})$
with one generator in degree $r$ for each $r$-cell of $X$, and
$$H_*(X;\Lambda)~=~H_*(\Lambda\otimes_{\Z[\pi_1(X)]}C(\widetilde{X}))~.$$

\subsection{Finite domination} \label{Finite domination}

A topological space $X$
is {\it finitely dominated} if there exist a finite $CW$ complex $K$,
maps $f:X \to K$, $g:K \to X$ and a homotopy $gf \simeq 1:X \to X$.
The finiteness obstruction of Wall \cite{Wall1} is a reduced projective
class $[X] \in \widetilde{K}_0(\Z[\pi_1(X)])$
such that $[X]=0$ if and only if $X$ is homotopy equivalent to a
finite $CW$ complex.

In the applications of the finiteness
obstruction to manifold topology $X=\overline{M}$ is an infinite
cyclic cover of a compact manifold $M$ -- see Chapter 17 of Hughes
and Ranicki \cite{HR} for the geometric wrapping up procedure which
shows that in dimension $\geqslant 5$ every tame manifold end has a
neighbourhood which is a finitely dominated infinite cyclic
cover $\overline{M}$ of a compact manifold $M$. Let $f:M \to S^1$
be a classifying map, so that $\overline{M}=f^*\R$, and let
$\overline{M}^+=f^*\R^+$. The finiteness
obstruction $[\overline{M}^+] \in \widetilde{K}_0(\Z[\pi_1(\overline{M})])$
is the end obstruction of Siebenmann \cite{Si}, such that
$[\overline{M}^+]=0$ if and only if the tame end can be closed, i.e.
compactified by a manifold with boundary.

Given a ring $A$ let $\Omega$ be the set of square matrices $\omega
\in M_r(A[z,z^{-1}])$ over the Laurent polynomial extension $A[z,z^{-1}]$
such that the $A$-module
$$P~=~{\rm coker}(\omega:A[z,z^{-1}]^r \to A[z,z^{-1}]^r)$$
is f.g. projective. The noncommutative Fredholm localization
$\Omega^{-1}A[z,z^{-1}]$ has
the universal property that a finite f.g. free $A[z,z^{-1}]$-module
chain complex $C$ is $A$-module chain equivalent to a finite f.g.
projective $A$-module chain complex if and only if $H_*(\Omega^{-1}C)=0$
(Ranicki \cite[Proposition 13.9]{R2}), with
$\Omega^{-1}C=\Omega^{-1}A[z,z^{-1}]\otimes_{A[z,z^{-1}]}C$.

Let $M$ be a connected finite $CW$ complex with a connected
infinite cyclic cover $\overline{M}$.
The fundamental group $\pi_1(M)$ fits into an extension
$$\{1\} \to \pi_1(\overline{M}) \to \pi_1(M) \to \Z \to \{1\}$$
and $\Z[\pi_1(M)]$ is a twisted Laurent polynomial extension
$$\Z[\pi_1(M)]~=~\Z[\pi_1(\overline{M})]_{\alpha}[z,z^{-1}]$$
with
$$\alpha~:~\pi_1(\overline{M}) \to \pi_1(\overline{M})~;~g \mapsto z^{-1}gz$$
the monodromy automorphism.  For the sake of simplicity 
only the untwisted case $\alpha=1$ will be considered here, so that
$\pi_1(M)=\pi_1(\overline{M}) \times \Z$.  The infinite cyclic cover
$\overline{M}$ is finitely dominated if and only if
$H_*(M;\Omega^{-1}\Z[\pi_1(M)])=0$, with
$A=\Z[\pi_1(\overline{M})]$ and $\Z[\pi_1(M)]=A[z,z^{-1}]$.  The
Farrell-Siebenmann obstruction $\Phi(M) \in {\rm Wh}(\pi_1(M))$ of an
$n$-dimensional manifold $M$ with finitely dominated infinite cyclic
cover $\overline{M}$ is such that $\Phi(M)=0$ if (and for $n \geqslant
6$ only if) $M$ is a fibre bundle over $S^1$ -- see \cite[Proposition
15.16]{R2} for the expression of $\Phi(M)$ in terms of the
$\Omega^{-1}\Z[\pi_1(M)]$-coefficient Reidemeister-Whitehead torsion
$$\tau(M;\Omega^{-1}\Z[\pi_1(M)])~=~\tau(\Omega^{-1}C(\widetilde{M}))
\in K_1(\Omega^{-1}\Z[\pi_1(M)])~.$$

\subsection{Codimension 1 splitting}

Surgery theory asks whether a homotopy equivalence of manifolds is
homotopic (or $h$-cobordant) to a homeomorphism -- in general, the
answer is no.  There are obstructions in the topological
$K$-theory of vector bundles, in the algebraic $K$-theory of
modules and in the algebraic $L$-theory of quad\-ratic forms.  The
algebraic $K$-theory obstruction lives in the Whitehead group
$Wh(\pi)$ of the fundamental group $\pi$.  The $L$-theory
obstruction lives in one of the surgery groups $L_*(\Z[\pi])$ of
Wall \cite{Wall2}, and is defined when the topological and
algebraic $K$-theory obstructions vanish.  The groups
$L_*(\Lambda)$ are defined for any ring with involution $\Lambda$
to be the generalized Witt groups of non\-singular quad\-ratic
forms over $\Lambda$.  For manifolds of dimension $\geqslant 5$
the vanishing of the algebraic obstructions is both a necessary
and sufficient condition for deforming a homotopy equivalence to a
homeomorphism.  See Ranicki \cite{R1} for the reduction of the
Browder-Novikov-Sullivan-Wall surgery theory to algebra.

A homotopy equivalence of $m$-dimensional manifolds $f:M' \to M$
{\it splits} along a submanifold $N^n \subset M^m$ if $f$ is
homotopic to a map (also denoted by $f$) such that $N'=f^{-1}(N)
\subset M'$ is also a submanifold, and the restriction $f\vert:N'
\to N$ is also a homotopy equivalence.  For codimension $m-n
\geqslant 3$ the splitting obstruction is just the ordinary
surgery obstruction $\sigma_*(f\vert) \in L_m(\Z[\pi_1(N)])$. For
codimension $m-n=1,2$ the splitting obstructions involve the
interplay of the knotting properties of codimension $(m-n)$
submanifolds and Mayer-Vietoris-type decompositions of the
algebraic $K$- and $L$-groups of $\Z[\pi_1(M)]$ in terms of the
groups of $\Z[\pi_1(N)]$, $\Z[\pi_1(M \backslash N)]$.

In the case $m-n=1$ $\pi_1(M)$ is a generalized free product, i.e. 
either an amalgamated free product or an $HNN$ extension, by the
Seifert-van Kampen theorem.  Codimension 1 splitting theorems and the
algebraic $K$- and $L$-theory of generalized free products are a major
ingredient of high-dimensional manifold topology, featuring in the work
of Stallings, Browder, Novikov, Wall, Siebenmann, Farrell, Hsiang,
Shaneson, Casson, Waldhausen, Cappell, $\dots$, and the author. 
Noncommutative localization provides a systematic development of this
algebra, using the intuition afforded by the topological applications
-- see Part 2 below for a more detailed discussion.

\subsection{Homology surgery theory}\label{Homology surgery theory}

For a morphism of
rings with involution ${\mathcal F}:\Z[\pi] \to \Lambda$ Cappell and
Shaneson \cite{CS1} considered the problem of whether a
$\Lambda$-coefficient homology equivalence of manifolds with
fundamental group $\pi$ is $H$-cobordant to a homeomorphism.  Again, the
answer is no in general, with obstructions in the topological
$K$-theory of vector bundles and in the homology surgery groups
$\Gamma_*({\mathcal F})$, which are generalized Witt groups of
$\Lambda$-non\-singular quad\-ratic forms over $\Z[\pi]$.
Vogel \cite{V1}, \cite{V2} identified the $\Lambda$-coefficient homology
surgery groups with the ordinary $L$-groups of the localization
$\Sigma^{-1}\Z[\pi]$ of $\Z[\pi]$ inverting the set $\Sigma$ of
$\Lambda$-invertible square matrices over $\Z[\pi]$
$$\Gamma_*({\mathcal F})~=~L_*({\Sigma}^{-1}\Z[\pi])~,$$
and identified the relative $L$-groups
$L_*(\Z[\pi] \to \Sigma^{-1}\Z[\pi])$
in the localization exact sequence
$$\dots \to L_n(\Z[\pi]) \to L_n({\Sigma}^{-1}\Z[\pi]) \to
L_n(\Z[\pi] \to {\Sigma}^{-1}\Z[\pi]) \to L_{n-1}(\Z[\pi]) \to \dots$$
with generalized Witt groups $L_*(\Z[\pi],\Sigma)$ of non\-singular
${\Sigma}^{-1}\Z[\pi]/\Z[\pi])$-valued quad\-ratic
linking forms on $\Sigma$-torsion $\Z[\pi]$-modules of homological
dimension 1.

\subsection{Codimension 2 embeddings}

Suppose given a codimension 2 embedding $N^n \subset M^{n+2}$ such
as a knot or link. Let $\Sigma^{-1}A$ be the localization of
$A=\Z[\pi_1(M \backslash N)]$ inverting the set $\Sigma$ of
matrices over $A$ which become invertible over $\Z[\pi_1(M)]$. By
Alexander duality the $\Sigma^{-1}A$-coefficient homology modules
$$H_*(M \backslash N;\Sigma^{-1}A)~\cong~H^{n+2-*}(M,N;\Sigma^{-1}A)
~~(* \neq 0,n+2)$$
are determined by the homotopy class of the inclusion $N \subset M$.
The $A$-coefficient homology groups $H_*(M \backslash N;A)$ and their
Poincar\'e duality properties reflect more subtle invariants of $N
\subset M$ such as knotting.  See Ranicki \cite{R2} for a general
account of high-dimensional codimension 2 embedding theory, including
some of the applications of noncommutative localization.

\subsection{Open books}  An $(n+2)$-dimensional manifold
$M^{n+2}$ is an {\it open book} if there exists a codimension 2
submanifold $N^n \subset M^{n+2}$ such that the complement
$M \backslash N$ is a fibre bundle over $S^1$.
Every odd-dimensional manifold is an open book. Quinn \cite{Q1} showed that
for $k \geqslant 2$ a $(2k+2)$-dimensional manifold $M$ is an open book if
and only if an asymmetric form over $\Z[\pi_1(M)]$ associated to $M$
represents 0 in the Witt group.  This obstruction was identified in
Ranicki \cite{R2} with an element in the $L$-group
$L_{2k+2}({\Omega}^{-1}\Z[\pi_1(M)][z,z^{-1}])$ of the Fredholm
localization of $\Z[\pi_1(M)][z,z^{-1}]$
(cf. section \ref{Finite domination} above).

\subsection{Boundary link cobordism}
An $n$-dimensional $\mu$-component boundary link is a codimension 2 embedding
$$N^n~=~\bigcup\limits_{\mu}S^n \subset M^{n+2}~=~S^{n+2}$$
with a $\mu$-component Seifert surface, in which case the
fundamental group of the complement $X=M\backslash N$ has
a compatible surjection $\pi_1(X) \to F_{\mu}$ onto
the free group on $\mu$ generators. Duval \cite{D} used the work of
Cappell and Shaneson \cite{CS2} and Vogel \cite{V2} to identify
the cobordism group of $n$-dimensional $\mu$-component boundary links
for $n \geqslant 2$ with the relative $L$-group
$L_{n+3}(\Z[F_{\mu}],\Sigma)$ in the
localization exact sequence
$$\dots \to L_{n+3}(\Z[F_{\mu}]) \to L_{n+3}(\Sigma^{-1}\Z[F_{\mu}])
\to L_{n+3}(\Z[F_{\mu}],\Sigma) \to L_{n+2}(\Z[F_{\mu}])\to \dots$$
with $\Sigma$ the set of $\Z$-invertible square matrices over
$\Z[F_{\mu}]$.  The even-dimensional boundary link cobordism groups are
$L_{2*+1}(\Z[F_{\mu}],\Sigma)=0$.  The cobordism class in
$L_{2k+2}(\Z[F_{\mu}],\Sigma)$ of a $(2k-1)$-dimensional
$\mu$-component boundary link $\cup_{\mu}S^{2k-1} \subset S^{2k+1}$ was
identified with the Witt class of a
$\Sigma^{-1}\Z[F_{\mu}]/\Z[F_{\mu}]$-valued non\-singular
$(-1)^{k+1}$-quad\-ratic linking form on $H_k(X;\Z[F_{\mu}])$,
generalizing the Blanchfield pairing on the homology of the infinite
cyclic cover of a knot.  The localization $\Sigma^{-1}\Z[F_{\mu}]$ was
identified by Dicks and Sontag \cite{DS} and Farber and Vogel \cite{FV} 
with a ring of rational functions in $\mu$ noncommuting variables.  
The high odd-dimensional boundary link cobordism groups $L_{2*+2}(\Z[F_{\mu}],\Sigma)$ 
have been computed by Sheiham \cite{Sh2}.

\subsection{Circle-valued Morse theory}
Novikov \cite{N} proposed the study of the critical points of Morse
functions $f:M \to S^1$ on compact manifolds $M$.  The `Novikov
complex' $C(M,f)$ over $\Z((z))=\Z[[z]][z^{-1}]$ has one generator for
each critical point of $f$, and the `Novikov homology'
$$H_*(C(M,f))~=~H_*(M;\Z((z)))$$
provides lower bounds on the number of critical points of Morse functions
in the homotopy class of $f$, generalizing the inequalities of the
classical Morse theory of real-valued functions $M \to \R$.
Suppose given a Morse function $f:M \to S^1$ with $\overline{M}=f^*\R$ 
such that $\pi_1(M)=\pi_1(\overline{M}) \times \Z$ (for the sake
of simplicity). Let $\Sigma$ be the set of square matrices over
$\Z[\pi_1(\overline{M})][z]$ which become invertible over
$\Z[\pi_1(\overline{M})]$ under the augmentation $z \mapsto 0$.
There is a natural morphism from the localization to the completion
$$\Sigma^{-1}\Z[\pi_1(M)] \to \widehat{\Z[\pi_1(M)]}~=~
\Z[\pi_1(\overline{M})][[z]][z^{-1}]$$
which is an injection if $\pi_1(M)$ is abelian or $F_{\mu}$ (Dicks and
Sontag \cite{DS}, Farber and Vogel \cite{FV}), but may not be an
injection in general (Sheiham \cite{Sh1}).  See Pajitnov \cite{P},
Farber and Ranicki \cite{FR}, Ranicki \cite{R3}, and Cornea and Ranicki
\cite{CR} for the construction and properties of Novikov complexes of
$f$ over $\widehat{\Z[\pi_1(M)]}$ and $\Sigma^{-1}\Z[\pi_1(M)]$. 
Naturally, noncommutative localization also features in the more
general Morse theory of closed 1-forms -- see Novikov \cite{N} and
Farber \cite{F}.

\subsection{3- and 4-dimensional manifolds}
See Garoufalidis and Kricker \cite{GK}, Quinn \cite{Q2} for
applications of noncommutative localization in the topology of 3-
and 4-dimensional manifolds.

\setcounter{section}{2}
\setcounter{subsection}{0}

\section*{Part 2. The algebraic K- and L-theory of generalized free products
via noncommutative localization}

A generalized free product of groups (or rings) is either an
amalgamated free product or an $HNN$ extension.  The expressions of
Schofield \cite{Sc} of generalized free products as noncommutative
localizations of triangular matrix rings combine with the localization
exact sequences of Neeman and Ranicki \cite{NR} to provide more
systematic proofs of the Mayer-Vietoris decompositions of Waldhausen
\cite{Wald} and Cappell \cite{C} of the algebraic $K$- and $L$-theory
of generalized free products.  The topological motivation for these
proofs comes from a noncommutative localization interpretation of
the Seifert-van Kampen and Mayer-Vietoris theorems.
If $(M,N\subseteq M)$ is a two-sided pair of connected $CW$ complexes
the fundamental group $\pi_1(M)$ is a generalized free product:
an amalgamated free product if $N$ separates $M$, and an $HNN$ extension
otherwise. The morphisms $\pi_1(N) \to \pi_1(M\backslash N)$
determine a triangular $k \times k$ matrix ring $A$ with universal
localization the full $k \times k$ matrix ring
$\Sigma^{-1}A=M_k(\Z[\pi_1(M)])$ ($k=3$ in the separating case,
$k=2$ in the non-separating case), such that
the corresponding presentations of the $\Z[\pi_1(M)]$-module chain
complex $C(\widetilde{M})$ of the universal cover $\widetilde{M}$ is
the assembly of an $A$-module chain complex constructed from
the chain complexes $C(\widetilde{N})$, $C(\widetilde{M\backslash N})$
of the universal covers $\widetilde{N}$, $\widetilde{M \backslash N}$
of $N$, $M \backslash N$.
The two cases will be considered separately, in sections \ref{HNN},
\ref{amalg}.

\subsection{The algebraic $K$-theory of a noncommutative localization}

Given an injective noncommutative localization $A \to \Sigma^{-1}A$ let
$H(A,\Sigma)$ be the exact category of homological dimension 1
$A$-modules $T$ which admit a f.g. projective $A$-module resolution
$$\xymatrix{0 \ar[r] & P \ar[r]^s & Q \ar[r] & T \ar[r] & 0}$$
such that $1 \otimes s : \Sigma^{-1}P \to \Sigma^{-1}Q$ is an
$\Sigma^{-1}A$-module isomorphism.
The algebraic $K$-theory localization exact sequence of
Schofield \cite[Theorem 4.12]{Sc}
$$K_1(A) \to K_1(\Sigma^{-1}A) \to K_1(A,\Sigma) \to K_0(A) \to K_0(\Sigma^{-1}A)$$
was obtained for any injective noncommutative localization $A \to \Sigma^{-1}A$,
with $K_1(A,\Sigma)\allowbreak =K_0(H(A,\Sigma))$.
Neeman and Ranicki \cite{NR} proved that if $A \to \Sigma^{-1}A$ is injective
and `stably flat'
$${\rm Tor}^A_i(\Sigma^{-1}A,\Sigma^{-1}A)~=~0~~(i \geqslant 1)$$
then
\begin{itemize}
\item[(i)] $\Sigma^{-1}A$ has the chain complex lifting property~:
every finite f.g. free $\Sigma^{-1}A$-module chain complex $C$ is chain equivalent
to $\Sigma^{-1}B$ for a finite f.g. projective $A$-module chain complex $B$,
\item[(ii)] the localization exact sequence extends to the higher $K$-groups
$$\dots \to K_n(A) \to K_n(\Sigma^{-1}A) \to K_n(A,\Sigma) \to K_{n-1}(A)
\to \dots \to K_0(\Sigma^{-1}A)$$
with $K_n(A,\Sigma)=K_{n-1}(H(A,\Sigma))$.
\end{itemize}

\subsection{Matrix rings}\label{triang}

The amalgamated free product of rings and the $HNN$ construction are
special cases of the following type of noncommutative localization
of triangular matrix rings.

Given rings $A_1,A_2$ and an $(A_1,A_2)$-bimodule $B$ define the triangular
$2\times 2$ matrix ring
$$A~=~\begin{pmatrix} A_1 & B \\ 0 & A_2 \end{pmatrix}~.$$
An $A$-module can be written as
$$M~=~\begin{pmatrix} M_1 \\ M_2 \end{pmatrix}$$
with $M_1$ an $A_1$-module, $M_2$ an $A_2$-module, together with an
$A_1$-module morphism $B\otimes_{A_2}M_2 \to M_1$.
The injection
$$A_1 \times A_2 \to A~;~(a_1,a_2) \mapsto
\begin{pmatrix} a_1 & 0 \\ 0 & a_2 \end{pmatrix}$$
induces isomorphisms of algebraic $K$-groups
$$K_*(A_1) \oplus K_*(A_2)~\cong~K_*(A)~.$$
The columns of $A$ are f.g. projective $A$-modules
$$P_1~=~\begin{pmatrix} A_1 \\ 0 \end{pmatrix}~,~P_2~=~
\begin{pmatrix} B \\ A_2 \end{pmatrix}$$
such that
$$\begin{array}{l}
P_1 \oplus P_2 ~ =~ A~,~{\rm Hom}_A(P_i,P_i)~=~A_i~(i=1,2)~,\\[1ex]
{\rm Hom}_A(P_1,P_2)~=~B~,~{\rm Hom}_A(P_2,P_1)~=~0~.
\end{array}$$
The noncommutative localization of $A$ inverting a non-empty
subset $\Sigma \subseteq {\rm Hom}_A(P_1,P_2)=B$
is the $2 \times 2$ matrix ring
$$\Sigma^{-1}A~=~M_2(C)~=~\begin{pmatrix} C & C \\ C & C \end{pmatrix}$$
with $C$ the endomorphism ring of the induced f.g. projective
$\Sigma^{-1}A$-module $\Sigma^{-1}P_1\cong \Sigma^{-1}P_2$.
The Morita equivalence
$$\{\Sigma^{-1}A\hbox{\rm -modules}\} \to \{C\hbox{\rm -modules}\}~;~
L \mapsto (C~C)\otimes_{\Sigma^{-1}A}L$$
induces isomorphisms in algebraic $K$-theory
$$K_*(M_2(C))~\cong~K_*(C)~.$$
The composite of the functor
$$\{A\hbox{\rm -modules}\} \to \{\Sigma^{-1}A\hbox{\rm -modules}\}~;~
M \mapsto {\Sigma}^{-1}M~=~\Sigma^{-1}A\otimes_AM$$
and the Morita equivalence is the {\it assembly} functor
$$\begin{array}{l}
\{A\hbox{\rm -modules}\} \to \{C\hbox{\rm -modules}\}~;\\[1ex]
M~=~\begin{pmatrix} M_1 \\ M_2 \end{pmatrix} \mapsto (C~C)\otimes_AM\\[1ex]
\hphantom{M~=~\begin{pmatrix} M_1 \\ M_2 \end{pmatrix} \mapsto}
=~{\rm coker}(C \otimes_{A_1}B \otimes_{A_2}M_2 \to C\otimes_{A_1} M_1
\oplus C\otimes_{A_2} M_2)
\end{array}$$
inducing the morphisms
$$K_*(A)~=~K_*(A_1)\oplus K_*(A_2) \to K_*(\Sigma^{-1}A)~=~K_*(C)$$
in the algebraic $K$-theory localization exact sequence.

There are evident generalizations to $k \times k$ matrix
rings for any $k \geqslant 2$.

\subsection{$HNN$ extensions}\label{HNN}
The $HNN$ extension $R*_{\alpha,\beta}\{z\}$ is defined for any
ring morphisms $\alpha,\beta:S \to R$, with
$$\alpha(s)z~=~z\beta(s) \in R*_{\alpha,\beta}\{z\}~~(s \in S)~.$$
Define the triangular $2 \times 2$ matrix ring
$$A~=~\begin{pmatrix} R & R_{\alpha} \oplus R_{\beta} \\ 0 & S \end{pmatrix}$$
with $R_{\alpha}$ the $(R,S)$-bimodule $R$ with $S$ acting on $R$ via $\alpha$,
and similarly for $R_{\beta}$. Let $\Sigma=\{\sigma_1,\sigma_2\}
\subset {\rm Hom}_A(P_1,P_2)$, with
$$\sigma_1~=~\begin{pmatrix} (1,0) \\ 0 \end{pmatrix}~,~
\sigma_2~=~\begin{pmatrix} (0,1) \\ 0 \end{pmatrix}~:~
P_1~=~\begin{pmatrix} R \\ 0 \end{pmatrix}
\to P_2 ~=~\begin{pmatrix} R_{\alpha} \oplus R_{\beta} \\ S \end{pmatrix}~.$$
The $A$-modules $P_1,P_2$ are f.g. projective since $P_1\oplus P_2=A$.
Theorem 13.1 of \cite{Sc} identifies
$$\Sigma^{-1}A~=~M_2(R*_{\alpha,\beta}\{z\})~.$$

\noindent{\it Example}
Let $(M,N \subseteq M)$ be a non-separating pair of connected
$CW$ complexes such that $N$ is two-sided in $M$
(i.e. has a neighbourhood $N \times [0,1] \subseteq M$)
with $M \backslash N=M_1$ connected
$$M~=~M_1 \cup_{N \times \{0,1\}} N \times [0,1]$$
$$
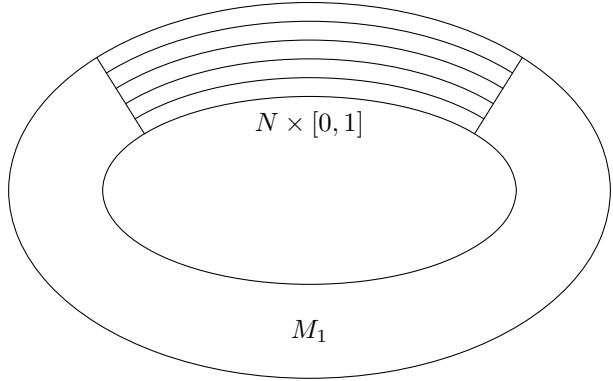
\epsfbox{b.ps}$$
By the Seifert-van Kampen theorem, the fundamental group $\pi_1(M)$ is
the $HNN$ extension determined by the morphisms
$\alpha,\beta:\pi_1(N) \to \pi_1(M_1)$ induced
by the inclusions $N\times \{0\} \to M_1$, $N\times \{1\} \to M_1$
$$\pi_1(M)~=~\pi_1(M_1)*_{\alpha,\beta}\{z\}~,$$
so that
$$\Z[\pi_1(M)]~=~\Z[\pi_1(M_1)]*_{\alpha,\beta}\{z\}~.$$
As above, define a triangular $2 \times 2$ matrix ring
$$A~=~\begin{pmatrix} \Z[\pi_1(N)] &
\Z[\pi_1(M_1)]_{\alpha} \oplus \Z[\pi_1(M_1)]_{\beta}\\
0 & \Z[\pi_1(M)] \end{pmatrix}$$
with noncommutative localization
$$\Sigma^{-1}A~=~M_2(\Z[\pi_1(M_1)]*_{\alpha,\beta}\{z\})~=~
M_2(\Z[\pi_1(M)])~.$$
Assume that $\pi_1(N) \to \pi_1(M)$ is injective,
so that the morphisms $\alpha,\beta$ are injective, and
the universal cover $\widetilde{M}$ is a union
$$\widetilde{M}~=~\bigcup\limits_{g\in [\pi_1(M):\pi_1(M_1)]}
g\widetilde{M}_1$$
of translates of the universal cover $\widetilde{M}_1$ of $M_1$, and
$$g_1\widetilde{M}_1 \cap g_2\widetilde{M}_1~=~
\begin{cases} h\widetilde{N}&\hbox{if $g_1\cap g_2z=h \in [\pi_1(M):\pi_1(N)]$}\\
g_1\widetilde{M}_1&\hbox{if $g_1=g_2$}\\
\emptyset&\hbox{if $g_1 \neq g_2$ and $g_1\cap g_2z=\emptyset$}
\end{cases}$$
with $h\widetilde{N}$ the translates of the universal cover
$\widetilde{N}$ of $N$.
In the diagram it is assumed that $\alpha,\beta$ are isomorphisms

\Draw
\LineAt(-160,-30,160,-30) \LineAt(-160,30,160,30)
\LineAt(-120,30,-120,-30) \LineAt(-40,30,-40,-30)
\LineAt(40,30,40,-30) \LineAt(120,30,120,-30) \MoveTo(-195,0)
\Text(--$\widetilde{M}$--) \MoveTo(-160,0) \Text(--$z^{-2}\widetilde{M}_1$--)
\MoveTo(-80,0) \Text(--$z^{-1}\widetilde{M}_1$--) \MoveTo(0,0)
\Text(--$\widetilde{M}_1$--) \MoveTo(80,0) \Text(--$z\widetilde{M}_1$--) \MoveTo(160,0)
\Text(--$z^2\widetilde{M}_1$--) \MoveTo(-120,-40) \Text(--$z^{-1}\widetilde{N}$--)
\MoveTo(-40,-40) \Text(--$\widetilde{N}$--) \MoveTo(40,-40) \Text(--$z\widetilde{N}$--)
    \MoveTo(120,-40) \Text(--$z^2\widetilde{N}$--) \EndDraw

\noindent
The cellular f.g. free chain complexes $C(\widetilde{M}_1)$, $C(\widetilde{N})$
are related by $\Z[\pi_1(M_1)]$-module chain maps
$$\begin{array}{l}
i_{\alpha}~:~\Z[\pi_1(M_1)]_{\alpha}\otimes_{\Z[\pi_1(N)]} C(\widetilde{N}) \to
C(\widetilde{M}_1)~,\\[1ex]
i_{\beta}~:~\Z[\pi_1(M_1)]_{\beta}\otimes_{\Z[\pi_1(N)]} C(\widetilde{N}) \to
C(\widetilde{M}_1)
\end{array}$$
defining a f.g. projective $A$-module chain complex
$\begin{pmatrix} C(\widetilde{M}_1) \\ C(\widetilde{N}) \end{pmatrix}$
with assembly the cellular f.g. free $\Z[\pi_1(M)]$-module chain complex of
$\widetilde{M}$
$$\begin{array}{l}
{\rm coker}\bigg(i_{\alpha}-zi_{\beta}~:\Z[\pi_1(M)]
\otimes_{\Z[\pi_1(N)]}C(\widetilde{N}) \to
\Z[\pi_1(M)]\otimes_{\Z[\pi_1(M_1)]}C(\widetilde{M}_1)\bigg)\\[2ex]
\hskip250pt =~C(\widetilde{M})
\end{array}$$
by the Mayer-Vietoris theorem.
\hfill$\qed$

Let $R*_{\alpha,\beta}\{z\}$ be an $HNN$ extension of rings in which
the morphisms $\alpha,\beta:S \to R$ are both injections of $(S,S)$-bimodule
direct summands, and $R_{\alpha},R_{\beta}$ are flat $S$-modules.
(This is the case in the above example if $\pi_1(N) \to \pi_1(M)$ is
injective). Then the natural ring morphisms
$$\begin{array}{l}
R \to R*_{\alpha,\beta}\{z\}~,~S \to R*_{\alpha,\beta}\{z\}~,\\[2ex]
A~=~\begin{pmatrix} R & R_{\alpha} \oplus R_{\beta} \\ 0 & S \end{pmatrix}
 \to \Sigma^{-1}A~=~M_2(R*_{\alpha,\beta}\{z\})
\end{array}$$
are injective, and $\Sigma^{-1}A$ is a stably flat universal
localization, with $H(A,\Sigma)= {\rm Nil}(R,S,\alpha,\beta)$ the
nilpotent category of Waldhausen \cite{Wald}.
The chain complex lifting property of $\Sigma^{-1}A$
gives a noncommutative localization proof of the existence of Mayer-Vietoris
presentations for finite f.g. free $R*_{\alpha,\beta}\{z\}$-module chain
complexes $C$
$$\xymatrix{0 \ar[r] & R*_{\alpha,\beta}\{z\}\otimes_SE
\ar[r]^-{i_{\alpha}-zi_{\beta}} & R*_{\alpha,\beta}\{z\}\otimes_RD
\ar[r] & C \ar[r] & 0}$$
with $D$ (resp. $E$) a finite f.g. free $R$- (resp. $S$-) module chain complex
(\cite{Wald}, Ranicki \cite{R4}).
The algebraic $K$-theory
localization exact sequence of \cite{NR}
$$\begin{array}{l}
\dots \to K_{n+1}(A,\Sigma)~=~
K_n(S) \oplus K_n(S) \oplus \widetilde{\rm Nil}_n(R,S,\alpha,\beta)\\
\xymatrix@C+15pt{
\ar[r]^-{\begin{pmatrix} \alpha &\beta & 0 \\ 1& 1 & 0\end{pmatrix}}&}
K_n(A)~=~K_n(R) \oplus K_n(S) \\[2ex]
\hskip150pt \to
 K_n(\Sigma^{-1}A)~=~K_n(R*_{\alpha,\beta}\{z\})\to \dots
\end{array}$$
is just the stabilization by $1:K_*(S) \to K_*(S)$ of the
Mayer-Vietoris exact sequence of \cite{Wald}
$$\xymatrix@C-5pt
{\dots \ar[r]& K_n(S) \oplus \widetilde{\rm Nil}_n(R,\alpha,\beta)}
\xymatrix@C+5pt{\ar[r]^-{(\alpha-\beta)\oplus 0}& K_n(R)}
\xymatrix@C-5pt
{\ar[r]& K_n(R*_{\alpha,\beta}\{z\})\ar[r]&\dots}$$

In particular, for $\alpha=\beta=1:S=R \to R$ the $HNN$ extension
is just the Laurent polynomial extension
$$R*_{\alpha,\beta}\{z\}~=~R[z,z^{-1}]$$
and the Mayer-Vietoris exact sequence splits to give the original splitting
of Bass, Heller and Swan \cite{BHS}
$$K_1(R[z,z^{-1}])~=~K_1(R) \oplus K_0(R) \oplus
\widetilde{\rm Nil}_0(R) \oplus\widetilde{\rm Nil}_0(R)$$
as well as its extension to the Quillen higher $K$-groups $K_*$.

\subsection{Amalgamated free products}\label{amalg}
The amalgamated free product $R_1*_SR_2$ is defined for any ring morphisms
$i_1:S \to R_1$, $i_2:S \to R_2$, with
$$r_1i_1(s)*r_2~=~r_1*i_2(s)r_2 \in R_1*_SR_2~~
(r_1 \in R_1,r_2\in R_2,s \in S)~.$$
Define the triangular $3 \times 3$ matrix ring
$$A~=~\begin{pmatrix} R_1 & 0 & R_1 \\
0 & R_2 & R_2 \\
0 & 0 & S \end{pmatrix}$$
and the $A$-module morphisms
$$\begin{array}{l}
\sigma_1~=~\begin{pmatrix} 1 \\ 0 \\ 0
\end{pmatrix}~:~P_1~=~\begin{pmatrix} R_1 \\ 0 \\ 0 \end{pmatrix} \to
P_3~=~\begin{pmatrix} R_1\\ R_2 \\ S \end{pmatrix}~,\\[4ex]
\sigma_2~=~\begin{pmatrix} 0 \\ 1 \\ 0
\end{pmatrix}~:~P_2~=~\begin{pmatrix} 0 \\ R_2 \\ 0 \end{pmatrix} \to
P_3~=~\begin{pmatrix} R_1\\ R_2 \\ S \end{pmatrix}~.
\end{array}$$
The $A$-modules $P_1,P_2,P_3$ are f.g. projective since
$P_1\oplus P_2 \oplus P_3=A$.
The noncommutative localization of $A$ inverting $\Sigma =
\{\sigma_1,\sigma_2\}$ is the full $3 \times 3$ matrix ring
$$\Sigma^{-1}A~=~M_3(R_1*_SR_2)$$
(a modification of Theorem 4.10 of \cite{Sc}).

\noindent{\it Example}
Let $(M,N \subseteq M)$ be a separating pair of $CW$ complexes such that $N$
has a neighbourhood $N \times [0,1] \subseteq M$ and
$$M~=~M_1 \cup_{N \times \{0\}} N\times [0,1]\cup_{N \times \{1\}}M_2$$
with $M_1,M_2,N$ connected.
$$
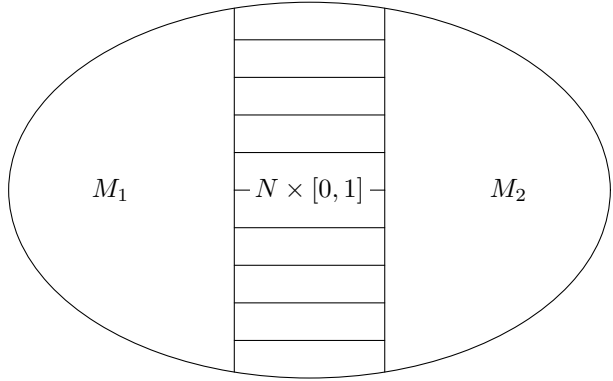
\epsfbox{a.ps}$$
By the Seifert-van Kampen theorem,
the fundamental group of $M$ is the  amalgamated free product
$$\pi_1(M)~=~\pi_1(M_1)*_{\pi_1(N)}\pi_1(M_2)~,$$
so that
$$\Z[\pi_1(M)]~=~\Z[\pi_1(M_1)]*_{\Z[\pi_1(N)]}\Z[\pi_1(M_2)]~.$$
As above, define a triangular matrix ring
$$A~=~\begin{pmatrix} \Z[\pi_1(M_1)] &  0 & \Z[\pi_1(M_1)]\\
0 & \Z[\pi_1(M_2)] & \Z[\pi_1(M_2)]  \\
0 & 0 & \Z[\pi_1(N)]\end{pmatrix}$$
with noncommutative localization
$$\Sigma^{-1}A~=~M_3(\Z[\pi_1(M_1)]*_{\Z[\pi_1(N)]}\Z[\pi_1(M_2)])~=~
M_3(\Z[\pi_1(M)])~.$$
Assume that $\pi_1(N) \to \pi_1(M)$ is injective, so that the morphisms
$$\begin{array}{l}
i_1~:~\pi_1(N) \to \pi_1(M_1)~,~i_2~:~\pi_1(N) \to \pi_1(M_2)~,\\[1ex]
\pi_1(M_1) \to \pi_1(M)~,~\pi_1(M_2) \to \pi_1(M)
\end{array}$$
are all injective, and the universal cover $\widetilde{M}$ of $M$ is a union
$$\widetilde{M}~=~\bigcup\limits_{g_1\in [\pi_1(M):\pi_1(M_1)]}
g_1\widetilde{M}_1
\cup_{\bigcup\limits_{h \in [\pi_1(M):\pi_1(N)]} h\widetilde{N}}
\bigcup\limits_{g_2\in [\pi_1(M):\pi_1(M_2)]}g_2\widetilde{M}_2$$
of $[\pi_1(M):\pi_1(M_1)]$ translates
of the universal cover $\widetilde{M}_1$ of $M_1$ and
$[\pi_1(M):\pi_1(M_2)]$ translates
of the universal cover $\widetilde{M}_2$ of $M_2$ with
intersection the $[\pi_1(M):\pi_1(N)]$ translates
of the universal cover $\widetilde{N}$ of $N$.
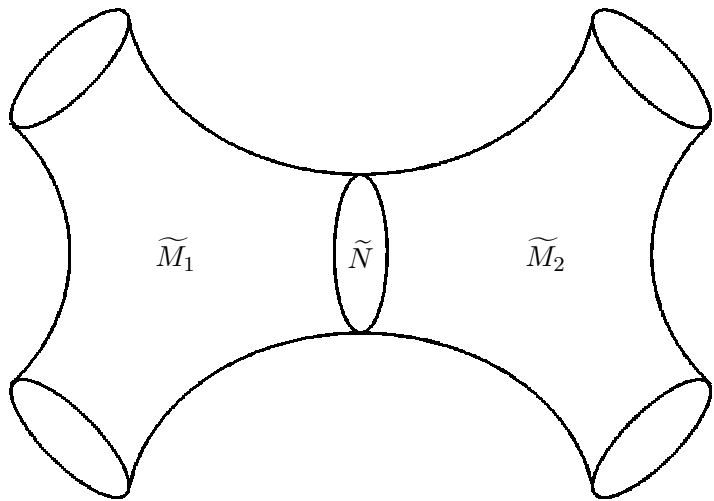
$$\Draw
\DrawOval(10,30)
\MoveTo(110,70)
\RotatedAxes(45,135)
\DrawOval(10,30)
\EndRotatedAxes
\MoveTo(-110,70)
\RotatedAxes(-45,-135)
\DrawOval(10,30)
\EndRotatedAxes
\MoveTo(110,-70)
\RotatedAxes(-45,-135)
\DrawOval(10,30)
\EndRotatedAxes
\MoveTo(-110,-70)
\RotatedAxes(45,135)
\DrawOval(10,30)
\EndRotatedAxes
\MoveTo(0,-100)
\DrawOvalArc(89,70)(8,172)
\MoveTo(0,100)
\DrawOvalArc(89,70)(188,352)
\MoveTo(187,0)
\DrawOvalArc(77,70)(138,222)
\MoveTo(-187,0)
\DrawOvalArc(77,70)(-42,42)
\MoveTo(0,0)
\Text(--$\widetilde{N}$--)
\MoveTo(70,0)
\Text(--$\widetilde{M}_2$--)
\MoveTo(-70,0)
\Text(--$\widetilde{M}_1$--)
\EndDraw$$

\noindent
The cellular f.g. free chain complexes $C(\widetilde{M}_1)$, $C(\widetilde{N})$ of the
universal covers $\widetilde{M}_1$, $\widetilde{N}$ are related by
$\Z[\pi_1(M_1)]$-module chain maps
$$\begin{array}{l}
i_1~:~\Z[\pi_1(M_1)]\otimes_{\Z[\pi_1(N)]} C(\widetilde{N}) \to
C(\widetilde{M}_1)~,\\[1ex]
i_2~:~\Z[\pi_1(M_2)]\otimes_{\Z[\pi_1(N)]} C(\widetilde{N}) \to
C(\widetilde{M}_2)
\end{array}$$
defining a f.g. projective $A$-module chain complex
$\begin{pmatrix} C(\widetilde{M}_1) \\ C(\widetilde{M}_2)
\\ C(\widetilde{N}) \end{pmatrix}$
with assembly the cellular f.g. free $\Z[\pi_1(M)]$-module chain complex of
$\widetilde{M}$
$$\begin{array}{l}
{\rm coker}\bigg(\begin{pmatrix} 1 \otimes i_1\\[1ex]
1 \otimes i_2 \end{pmatrix}:\Z[\pi_1(M)]\otimes_{\Z[\pi_1(N)]}C(\widetilde{N}) \to\\[1ex]
\hskip40pt
\Z[\pi_1(M)]\otimes_{\Z[\pi_1(M_1)]}C(\widetilde{M}_1)\oplus
\Z[\pi_1(M)]\otimes_{\Z[\pi_1(M_1)]}C(\widetilde{M}_2)\bigg)\\[2ex]
\hskip200pt =~ C(\widetilde{M})
\end{array}$$
by the Mayer-Vietoris theorem.\hfill$\qed$

Let $R_1*_SR_2$ be an amalgamated free product of rings in which
the morphisms $i_1:S \to R_1$, $i_2:S \to R_2$ are both injections of
$(S,S)$-bimodule direct summands, and $R_1,R_2$ are flat $S$-modules.
(This is the case in the above example if $\pi_1(N) \to \pi_1(M)$ is
injective).
Then the natural ring morphisms
$$\begin{array}{l}
R_1 \to R_1*_SR_2~,~R_2 \to R_1*_SR_2~,~S \to R_1*_SR_2~,\\[2ex]
A~=~\begin{pmatrix} R_1& 0 & R_1 \\
0 & R_2 & R_2 \\
0 & 0 & S \end{pmatrix} \to \Sigma^{-1}A~=~M_3(R_1*_SR_2)
\end{array}$$
are injective, and $\Sigma^{-1}A$ is a stably flat
noncommutative localization, with $H(A,\Sigma)=
{\rm Nil}(R_1,R_2,S)$ the nilpotent category of Waldhausen \cite{Wald}.
The chain complex lifting property of $\Sigma^{-1}A$
gives a noncommutative localization proof of the existence of Mayer-Vietoris
presentations for finite f.g. free $R_1*_SR_2$-module chain
complexes $C$
$$\xymatrix@C-6pt{0 \ar[r] & R_1*_SR_2\otimes_SE
\ar[r] & R_1*_SR_2 \otimes_{R_1}D_1 \oplus R_1*_SR_2 \otimes_{R_2}D_2
\ar[r] & C \ar[r] & 0}$$
with $D_i$ (resp. $E$) a finite f.g. free $R_i$- (resp. $S$-) module chain complex
(\cite{Wald}, Ranicki \cite{R4}).
The algebraic $K$-theory localization exact sequence of \cite{NR}
$$\begin{array}{l}
\dots \to K_{n+1}(A,\Sigma)~=~
K_n(S) \oplus K_n(S) \oplus \widetilde{\rm Nil}_n(R_1,R_2,S)\\
\xymatrix@C+15pt{
\ar[r]^-{\begin{pmatrix} i_1 & 0 & 0 \\ 0& i_2 & 0 \\
1 & 1 & 0 \end{pmatrix}}&}
K_n(A)~=~K_n(R_1) \oplus K_n(R_2) \oplus K_n(S)\\[2ex]
\hskip150pt \to  K_n(\Sigma^{-1}A)~=~K_n(R_1*_SR_2)\to \dots
\end{array}$$
is just the stabilization by $1:K_*(S) \to K_*(S)$ of the
Mayer-Vietoris exact sequence of \cite{Wald}
$$\begin{array}{l}
\xymatrix@C-5pt
{\dots \ar[r]& K_n(S) \oplus \widetilde{\rm Nil}_n(R_1,R_2,S)}\\
\xymatrix@C+15pt{
\ar[r]^-{\begin{pmatrix} i_1 & 0 \\ i_2 & 0 \end{pmatrix}}
&K_n(R_1) \oplus K_n(R_2) \ar[r]& K_n(R_1*_SR_2)\ar[r]&\dots}
\end{array}$$

\subsection{The algebraic $L$-theory of a noncommutative localization}

See Chapter 3 of Ranicki \cite{R0} for the algebraic $L$-theory
of a commutative localization.

The algebraic $L$-theory of a ring $A$ depends on an involution,
that is a function $\raise4pt\hbox{$\overline{~}$}:A \to A;a \mapsto \overline{a}$
such that
$$\overline{a+b}~=~\overline{a}+\overline{b}~,~
\overline{ab}~=~\overline{b}\,\overline{a}~,~\overline{\overline{a}}~=~a~,~
\overline{1}~=~1~~(a,b\in A)~.$$
For an injective noncommutative localization $A \to \Sigma^{-1}A$ of
a ring $A$ with an involution which extends to $\Sigma^{-1}A$
Vogel \cite{V2} obtained a localization exact sequence in quadratic $L$-theory
$$\dots \to L_n(A) \to L_n(\Sigma^{-1}A) \to L_n(A,\Sigma) \to L_{n-1}(A)
\to \dots$$
with $L_n(A,\Sigma)=L_{n-1}(H(A,\Sigma))$. (See \cite{NR} for the
symmetric $L$-theory localization exact sequence in the stably flat case).
At first sight, it does not appear possible to apply this sequence
to the triangular matrix rings of sections \ref{triang}, \ref{HNN},
\ref{amalg}. How does one define an involution on a triangular
matrix ring
$$A~=~\begin{pmatrix} A_1 & B\\ 0 & A_2 \end{pmatrix}~?$$
The trick is to observe that if $A_1,A_2$ are rings with involution,
and $(B,\beta)$ is a nonsingular symmetric form over $A_1$ such that
$B$ is an $(A_1,A_2)$-bimodule then $A$ has a {\it chain duality}
in the sense of Definition 1.1 of Ranicki \cite{R1}, sending an
$A$-module $M=\begin{pmatrix} M_1 \\ M_2 \end{pmatrix}$ to the
1-dimensional $A$-module chain complex
$$TM~:~TM_1~=~\begin{pmatrix} M^*_1 \\ 0 \end{pmatrix} \to
TM_0~=~\begin{pmatrix} B\otimes_{A_2}M^*_2 \\ M_2^* \end{pmatrix}~.$$
The quadratic $L$-groups of $A$ are just the relative $L$-groups
in the exact sequence
$$\dots \to L_n(A) \to L_n(A_2)
\xymatrix{\ar[r]^{(B,\beta)\otimes_{A_2}-}&} L_n(A_1) \to L_{n-1}(A)
\to \dots~.$$
In particular, for generalized free products of rings with
involution the triangular matrix rings $A$ of section \ref{HNN}, \ref{amalg}
have such chain dualities, and in the injective case the torsion
$L$-groups $L_*(A,\Sigma)=L_{*-1}(H(A,\Sigma))$
in the localization exact sequence
$$\dots \to L_n(A) \to L_n(\Sigma^{-1}A) \to L_n(A,\Sigma) \to L_{n-1}(A)
\to \dots$$
are just the unitary nilpotent $L$-groups UNil$_*$ of
Cappell \cite{C}.

\noindent School of Mathematics\newline
University of Edinburgh \newline
James Clerk Maxwell Building \newline
King's Buildings \newline
Mayfield Road \newline
Edinburgh EH9 3JZ\newline
SCOTLAND, UK
\bigskip

\noindent e-mail aar@maths.ed.ac.uk

\end{document}